\input epsf
\documentstyle{amsppt}
\pagewidth{6.4truein}\hcorrection{0in}
\pageheight{9truein}\vcorrection{0.00in}

\TagsOnRight
\NoRunningHeads
\catcode`\@=11
\def\logo@{}
\footline={\ifnum\pageno>1 \hfil\folio\hfil\else\hfil\fi}
\topmatter
\title The interaction of collinear gaps of arbitrary charge in a two dimensional dimer system
\endtitle
\author Mihai Ciucu\endauthor
\thanks Research supported in part by NSF grant DMS-0801625.
\endthanks
\affil
  Department of Mathematics, Indiana University\\
  Bloomington, Indiana 47405-5701
\endaffil
\abstract 
The correlation of gaps in dimer systems was introduced in 1963 by Fisher and Stephenson, who looked at the interaction of two monomers generated by the rigid exclusion of dimers on the closely packed square lattice. In previous work we considered the analogous problem on the hexagonal lattice, and we extended the set-up to include the correlation of any finite number of monomer clusters. For fairly general classes of monomer clusters we proved that the asymptotics of their correlation is given, for large separations between the clusters, by a multiplicative version of Coulomb's law for 2D electrostatics. However, our previous results required that the monomer clusters consist (with possibly one exception) of an even number of monomers. In this paper we determine the asymptotics of general defect clusters along a lattice diagonal in the square lattice (involving an arbitrary, even or odd number of monomers), and find that it is given by the same Coulomb law. We also obtain a conceptual interpretation for the multiplicative constant as the product of the correlations of the individual clusters.

\endabstract 
\endtopmatter
\document

\def\mysec#1{\bigskip\centerline{\bf #1}\message{ * }\nopagebreak\par\bigskip}

\def\myref#1{\item"{[{\bf #1}]}"} 
 
\def\pf{{\it Proof.\ }} 

\def\epf{\hfill{$\square$}\smallpagebreak}

\def\cite#1{\relaxnext@
  \def\nextiii@##1,##2\end@{[{\bf##1},\,##2]}%
  \in@,{#1}\ifin@\def\next{\nextiii@#1\end@}\else
  \def\next{[{\bf#1}]}\fi\next}
\def\proclaimheadfont@{\smc}

\def\pf{{\it Proof.\ }}

\define\Z{{\Bbb Z}}

\define\M{\operatorname{M}}

\define\q{\operatorname{q}}

\define\com{\operatorname{c}}

\define\de{\operatorname{d}}

\define\co{\operatorname{c}}
\define\ho{\operatorname{h}}
\define\se{\operatorname{s}}
\define\st{\operatorname{st}}
\define\De{\operatorname{D}}
\define\s{\operatorname{supp}}
\define\te{\operatorname{t}}
\define\twoline#1#2{\line{\hfill{\smc #1}\hfill{\smc #2}\hfill}}
\define\threeline#1#2#3{\line{\hfill{\smc #1}\hfill{\smc #2}\hfill{\smc #3}\hfill}}
\define\ltwoline#1#2{\line{{\smc #1}{\smc #2}}}

\def\mypic#1{\epsffile{figs/#1}}



  \define\Baikone{1}
  \define\Baiktwo{2}
  \define\FT{3}
  \define\ri{4} 
  \define\sc{5}
  \define\ppone{6}
  \define\ec{7}
  \define\ef{8}
  \define\ov{9}
  \define\CEP{10}
  \define\FS{11}
  \define\FisherIsing{12}
  \define\Glaish{13}
  \define\Hart{14}
  \define\SM{15}
  \define\MRR{16}
  \define\ZI{17}

%


\define\eac{1.1}
\define\ead{1.2}
\define\eae{1.3}
\define\eaf{1.4}
\define\eah{1.5}
\define\eai{1.6}
\define\eaj{1.7}

\define\ebba{2.1}
\define\ebbb{2.2}
\define\ebbc{2.3}
\define\ebaa{2.4}
\define\ebab{2.5}
\define\ebac{2.6}
\define\ebad{2.7}
\define\ebae{2.8}

\define\ebb{2.9}
\define\ebc{2.10}
\define\ebd{2.11}
\define\ebe{2.12}
\define\ebff{2.13}

\define\eca{3.1}
\define\ecb{3.2}

\define\ebh{4.1}
\define\ebf{4.2}
\define\ebj{4.3}
\define\ebk{4.4}
\define\ebkk{4.5}

\define\ebl{5.1}
\define\ebm{5.2}
\define\ebn{5.3}
\define\ebo{5.4}

\define\eea{5.5}
\define\ebnp{5.6}
\define\ebop{5.7}

\define\eboo{5.8}
\define\eboa{5.9}
\define\ebp{5.10}
\define\ebpa{5.11}
\define\ebpaa{5.12}
\define\ebpaaa{5.13}
\define\ebq{5.14}
\define\ebr{5.15}
\define\ebrr{5.16}
\define\ebs{5.17}
\define\ebt{5.18}
\define\ebuone{5.19}
\define\ebuu{5.20}
\define\ebu{5.21}
\define\ebv{5.22}
\define\ebw{5.23}
\define\ebx{5.24}

\define\efa{6.1}
\define\efb{6.2}
\define\efc{6.3}
\define\efd{6.4}
\define\efe{6.5}
\define\eff{6.6}
\define\efg{6.7}
\define\efh{6.8}
\define\efi{6.9}
\define\efj{6.10}
\define\efk{6.11}
\define\efii{6.12}
\define\efjj{6.13}
\define\efkk{6.14}
\define\efl{6.15}
\define\efm{6.16}
\define\efn{6.17}
\define\efo{6.18}
\define\efp{6.19}
\define\efq{6.20}
\define\efr{6.21}
\define\efs{6.22}
\define\eft{6.23}
\define\efu{6.24}
\define\efw{6.25}
\define\efmm{6.26}
\define\efoo{6.27}
\define\efpp{6.28}
\define\efqq{6.29}

\define\ega{7.1}
\define\egb{7.2}
\define\egc{7.3}
\define\egd{7.4}
\define\ege{7.5}
\define\egf{7.6}
\define\egg{7.7}
\define\egh{7.8}
\define\egi{7.9}
\define\egj{7.10}
\define\egk{7.11}
\define\egl{7.12}

\define\eha{8.1}
\define\ehb{8.2}
\define\ehc{8.3}
\define\ehd{8.4}
\define\ehe{8.5}
\define\ehf{8.6}
\define\ehg{8.7}
\define\ehh{8.8}
\define\ehi{8.9}
\define\ehj{8.10}
\define\ehk{8.11}
\define\ehl{8.12}
\define\ehm{8.13}
\define\ehn{8.14}
\define\eho{8.15}
\define\ehp{8.16}
\define\ehq{8.17}
\define\ehr{8.18}
\define\ehs{8.19}
\define\eht{8.20}
\define\ehu{8.21}
\define\ehv{8.22}
\define\ehw{8.23}

\define\taa{1.1}

\define\tba{2.1}

\define\tca{3.1}

\define\tbb{4.1}

\define\tbc{5.1}
\define\tbcp{5.1'}

\define\tbd{5.2}

\define\tfa{6.1}
\define\tfb{6.2}
\define\tfc{6.3}

\define\tfe{6.4}

\define\tga{7.1}

\define\faa{1.1}
\define\fab{1.2}
\define\fac{1.3}
\define\fad{1.4}
\define\fae{1.5}

\define\fba{5.1}
\define\fbb{5.2}
\define\fbc{5.3}
\define\fbd{5.4}

\define\fha{8.1}

\vskip-0.05in
\mysec{Introduction}

The correlation of gaps in dimer systems was introduced in 1963 by Fisher and Stephenson \cite{\FS}, who looked at the interaction of two monomers generated by the rigid exclusion of dimers completely covering the rest of the square lattice. In previous work (see \cite{\ri}, \cite{\sc}, \cite{\ec}, \cite{\ef} and \cite{\ov}) we considered the analogous problem on the hexagonal lattice, and we extended the set-up to include the correlation of any finite number of monomer clusters. For fairly general classes of monomer clusters we proved that the asymptotics of their correlation is given, for large separations between the clusters, by a multiplicative version of Coulomb's law for 2D electrostatics. However, our previous results required that the monomer clusters consist (with possibly one exception, see \cite{\sc}) of an even number of monomers. 

Odd monomer clusters (or odd gaps for short) are notoriously hard to handle. The simple case of two clusters, each consisting of a single monomer, is the still open conjecture of the rotational invariance of the monomer-monomer correlation, phrased in 1963 by Fisher and Stephenson. The only proved cases in the literature involving an odd gap are Hartwig's result \cite{\Hart} (on two monomers, one adjacent to the diagonal through the other) and our earlier result \cite{\sc} on a symmetric distribution of gaps on the hexagonal lattice (where a single monomer was allowed on the symmetry axis).

In the main result of this paper (see Theorem {\tba}) we determine the asymptotics of the correlation of general defect clusters along a lattice diagonal in the square lattice (involving an {\it arbitrary}, even or odd number of monomers), and find that it is given by the same Coulomb law. In order for our arguments to work, we modify Fisher and Stephenson's definition of the correlation, which they made by including the defects around the center of large squares, by including them around the center of large Aztec diamonds; this approach leads to the simple product formula of Theorem {\taa}, which is the starting point of our analysis (another key point is described in Remark 7, at the end of Section 6). However, based on several examples we worked out, we conjecture that our modified definition leads to precisely the same correlation values as Fisher and Stephenson's original definition (this is also in agreement with the results of Cohn, Elkies and Propp \cite{CEP} according to which, in the scaling limit, the domino statistics is undistorted at the very center of the Aztec diamond).
This proves in particular a conjecture of physicists Krauth and Moessner \cite{\SM}, who predicted, based on Monte-Carlo simulations, that two monomers of the same color in a dimer system on $\Z^2$ interact according to a Coulomb repulsion. An interesting additional feature of our result is that we also obtain a conceptual interpretation for the multiplicative constant as the product of the correlations of the individual clusters.

When specialized to the case when each defect cluster is a single monomer, our result can be viewed as a counterpart of a result of Zuber and Itzykson \cite{\ZI} determining the $n$-point spin correlation in the Ising model on the square lattice along a lattice line. Indeed, using Fisher's mapping \cite{\FisherIsing} of the Ising model on the dimer model, the latter is seen to be equivalent to a certain average of the correlation of collinear monomers in the lattice obtained by tiling the plane with equilateral triangles and regular dodecagons. 

Another point of contact with the previous literature is the detailed study of the correlation of collinear edges in large hexagonal regions on the hexagonal lattice, carried out by Baik, Kriecherbauer, McLaughlin and Miller in \cite{\Baikone} and \cite{\Baiktwo}.













\mysec{1. An exact enumeration result for Aztec rectangles with defects}

%
%
%
%
%
%
%

The results of \cite{\sc} and \cite{\ec} on the asymptotics of the correlation of holes in lozenge tilings were built upon exact product formulas giving the number of lozenge tilings of certain hexagonal regions enclosing the holes (see \cite{\ppone}). We provide in this section a counterpart of the latter formulas for the square lattice.

\topinsert
\twoline{\mypic{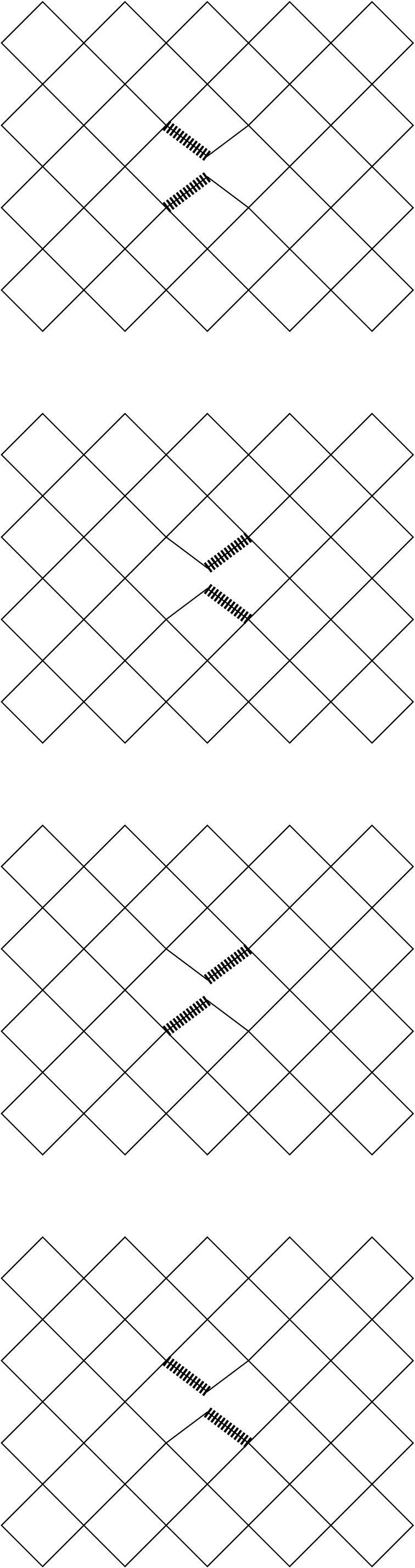}}{\mypic{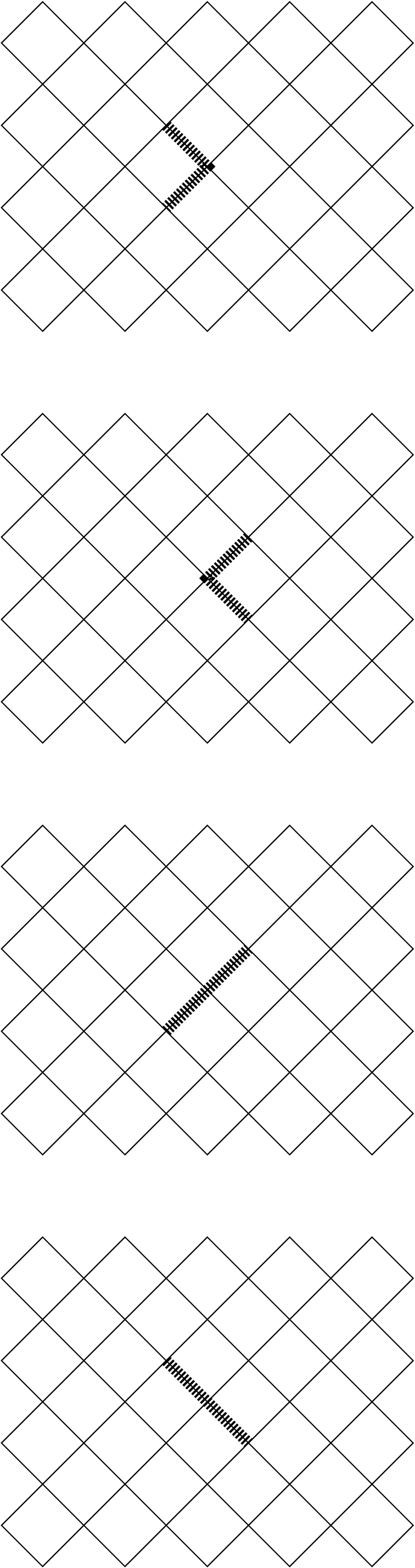}}
\medskip
\twoline{Figure {\faa}. {\rm (a) The four types of local matchings at a separation.}}
{{\rm \ (b) The corresponding trimers on $\Z^2$.}}
\endinsert

In addition to unit holes (or monomers), we consider also the following new type of defect on the square lattice. If $\Z^2$ is drawn so that the lattice lines form angles of $\pm\frac{\pi}{4}$ with the horizontal, we say that there is a {\it separation} at $v\in\Z^2$ if the vertex $v$ and its four incident edges are replaced by two new vertices, $v'$ and $v''$, and four new edges, two connecting $v'$ to the northern neighbors of $v$ and the other two connecting $v''$ to the southern neighbors of $v$ (see Figure {\faa(a)}).

Note that from the point of view of perfect matchings, a separation at $v\in\Z^2$ is equivalent to a superposition of four trimers, each centered at $v$ (see Figure {\faa}, and also the remark after Theorem {\tba}): two ``hooks'' (one pointing left, the other right) and two ``bars'' (one of slope 1, the other of slope $-1$). Indeed, $v'$ must be matched to one of the two northern neighbors of $v$, and $v''$ to one of the two southern neighbors of $v$; it is apparent that the four resulting combinations correspond to the four trimers described above.

Consider a $(2m+1)\times(2n+1)$ rectangular chessboard and suppose the corners are black. 
The {\it Aztec rectangle} $AR_{m,n}$ is the graph whose vertices are the white squares and whose edges connect precisely those pairs of white squares that are diagonally adjacent. 

Let $k,l \geq 0$ be integers, and consider the Aztec rectangle $AR_{2n,2n+k-l}$. Let $\ell$ be its horizontal symmetry axis, and label the vertices on $\ell$ from left to right by $1,2,\dotsc,2n+k-l$. Set $[m]:=\{1,2,\dotsc,m\}$. For any disjoint subsets $H,S\subseteq [2n+k-l]$ so that $|H|=k$ and $|S|=l$, define $AR_{2n,2n+k-l}(H,S)$ to be the graph obtained from $AR_{2n,2n+k-l}$ by deleting the vertices on $\ell$ whose labels belong to $H$, and creating separations at those vertices on $\ell$ whose labels belong to $S$ (see Figure {\fab} for an example). Clearly, $AR_{2n,2n+k-l}(H,S)$ is bipartite (i.e., its vertices can be colored black and white so that each edge has oppositely colored endpoints), and one readily checks that it is also balanced (i.e., the two color classes have the same number vertices), an obvious necessary condition for it to have perfect matchings.

\topinsert
\centerline{\mypic{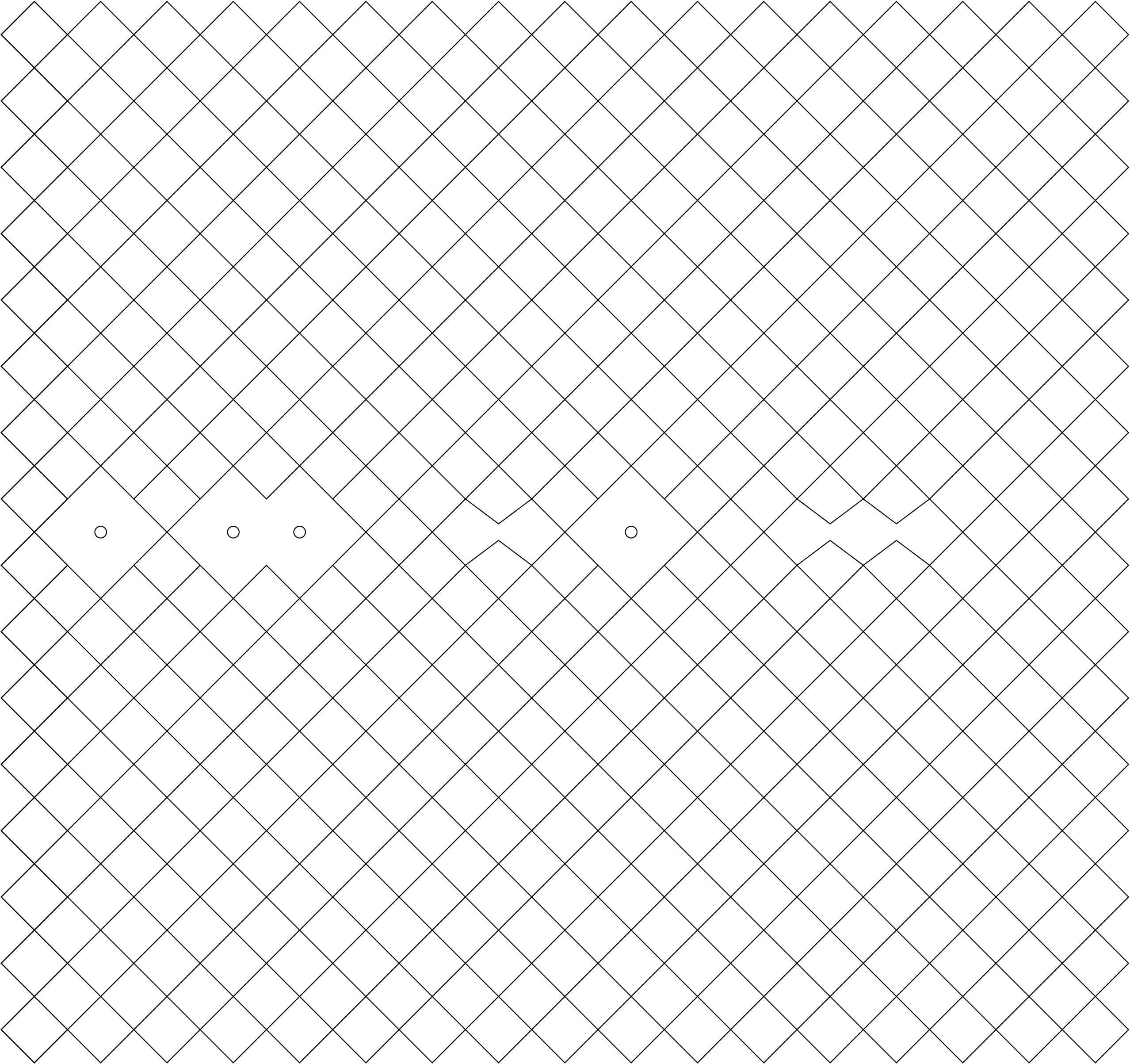}}
\centerline{{\smc Figure~{\fab}.} $AR_{16,17}(\{2,4,5,10\},\{8,13,14\})$.}
\endinsert

\topinsert
\centerline{\mypic{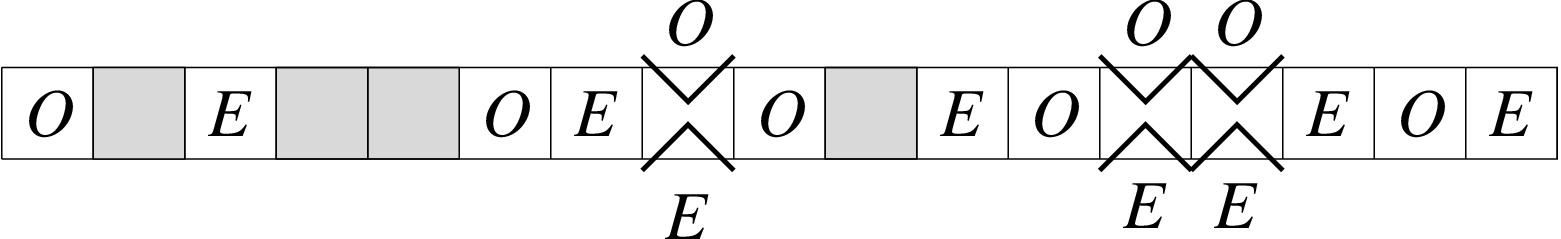}}
\centerline{{\smc Figure~{\fac}.} The $OE$-labeled strings corresponding to Figure {\fab}.}
\endinsert

Denote by $\M(G)$ the number of perfect matchings of the graph $G$.
The central exact product formula mentioned in the first paragraph of this section is the following.

\proclaim{Theorem \taa} For any integers $k,l \geq 0$ and disjoint subsets $H,S\subseteq [2n+k-l]$ with $|H|=k$ and $|S|=l$, the number of perfect matchings of $AR_{2n,2n+k-l}(H,S)$ is obtained as follows. Consider a string of $2n+k-l$ cells labeled from left to right $1,2,\dotsc,2n+k-l$. Focus on the $2n-2l$ cells whose labels are not in $H\cup S$, and write alternately $O$'s and $E$'s in them, starting from the left. Finally write both an $O$ and an $E$ in each cell whose label belongs to $S$ (Figure~{\fac} shows the labeling corresponding to the example in Figure~{\fab}). Let ${\Cal O}$ be the set of labels of the cells containing $O$, and ${\Cal E}$ the set of labels of the cells containing $E$. Then 
$$
\M(AR_{2n,2n+k-l}(H,S))=\frac{2^{n^2+2n-l}}{(0!\,1!\cdots(n-1)!)^2}\Delta({\Cal O})\Delta({\Cal E}),\tag\eac
$$
where for $T=\{t_1,\dotsc,t_m\}$, $t_1<\cdots<t_m$,
$$
\Delta(T):=\prod_{1\leq i<j\leq m}(t_j-t_i).\tag\ead
$$
\endproclaim

\flushpar
{\smc Remark 1.} The $O$-$E$-labeling of the string of $2n+k-l$ cells reflects the defects on the symmetry axis
of the Aztec rectangle in the following suggestive way (indicated in Figure~{\fad} for the labeling of Figure~{\fac}). If there are no monomers or separations, a string $OEOE\cdots OE$ of length $2n$ is obtained. The presence of each monomer introduces a unit {\it gap} in this string (cut the string and shift out one unit), while each separation introduces a unit {\it overlap} (cut the string and shift in one unit). The regular pattern $OEOE\cdots OE$ can be restored by closing in the gaps 
and undoing the overlaps.

\topinsert
\centerline{\mypic{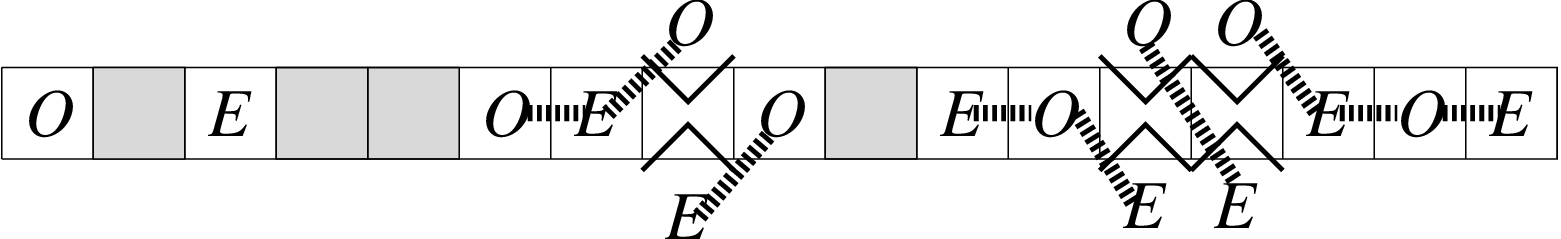}}
\centerline{{\smc Figure~{\fad}.} Holes and separations correspond to gaps and overlaps in an alternating $O$-$E$ string of length $2n$.}
\endinsert

\medskip
{\it Proof of Theorem 1.1.} Given a set $T=\{t_1,t_2,\dotsc,t_m\}$, $t_1<t_2<\cdots < t_m$, with $m$ even, define the subsets $T_o$ and $T_e$ by
$$
\spreadlines{2\jot}
\align
T_o:&=\{t_1,t_3,t_5,\dotsc,t_{m-1}\},\tag\eae
\\
T_e:&=\{t_2,t_4,t_6,\dotsc,t_m\}.\tag\eaf
\endalign
$$
For a set $T=\{t_1,\dotsc,t_m\}$ with $1\leq t_1<\cdots<t_m\leq n$ integers, let $\bar{R}_{m,n}(T)$ be the graph obtained from the Aztec rectangle $AR_{m,n}$ by deleting all its bottom vertices except for the $t_1$-th, $t_2$-th, ..., $t_m$-th (an example is illustrated in Figure {\fae}). As shown in \cite{\FT}, combinatorial arguments and \cite{\MRR, Theorem\,2} imply that the number of perfect matchings of this graph is given by
$$
\M(\bar{R}_{m,n}(T))=\frac{2^{m(m+1)/2}}{0!\,1!\,\cdots\,(m-1)!}\prod_{1\leq i<j\leq m}(t_j-t_i).\tag\eah
$$

\topinsert
\vskip0.2in
\centerline{\mypic{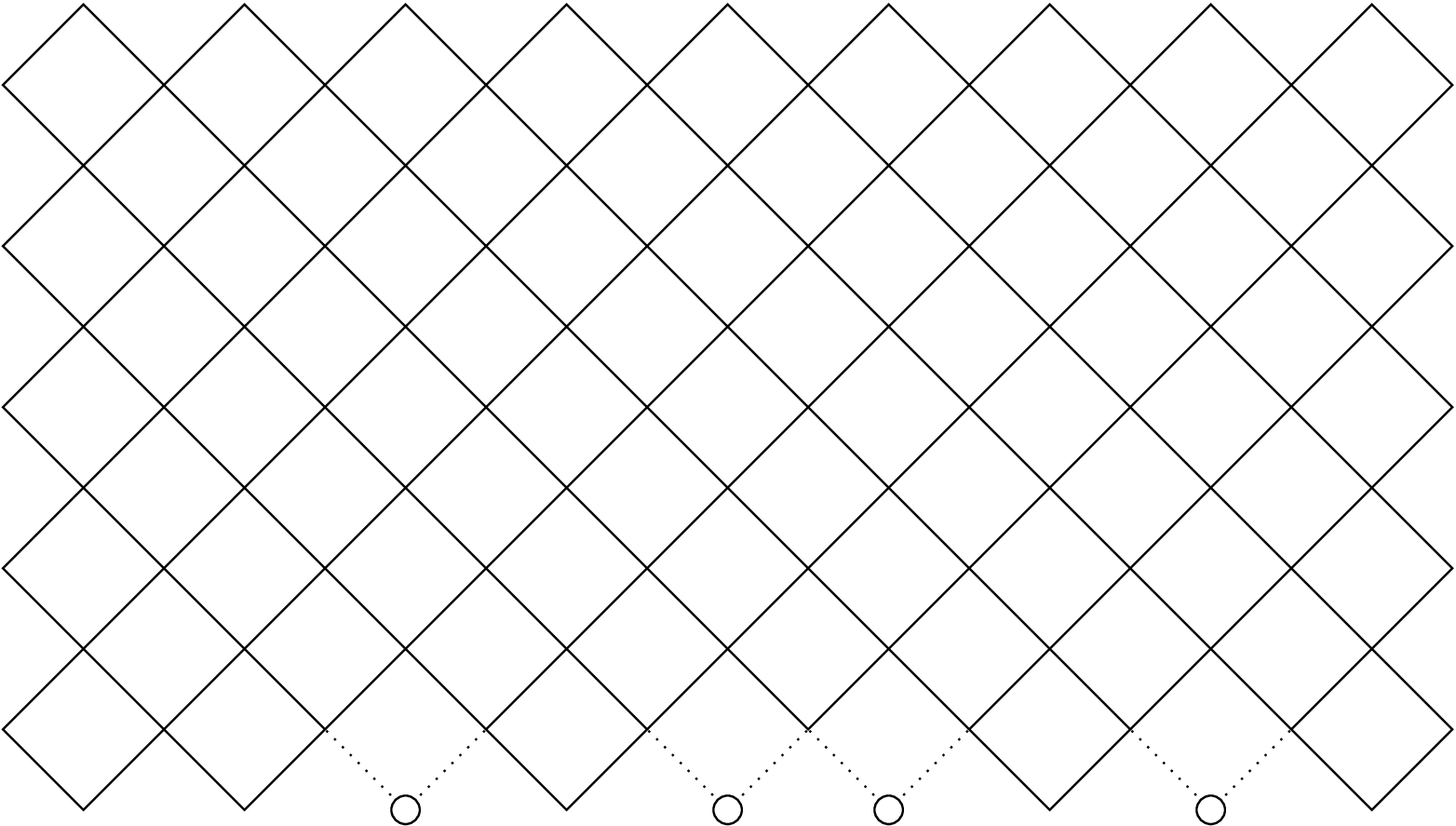}}
\medskip
\centerline{{\smc Figure~{\fae}.} $\bar{R}_{5,9}(\{1,2,4,7,9\})$ (deleted vertices of $AR_{5,9}$ are circled, deleted edges dotted).}
\endinsert

The graph $AR_{2n,2n+k-l}(H,S)$ is a planar bipartite graph that is symmetric with
respect to $\ell$. Thus the factorization theorem \cite{\FT,Theorem\,1.2} can be applied to it. This yields
$$
\M(AR_{2n,2n+k-l}(H,S))=2^{n-l}\M(\bar{R}_{n,2n+k-l}([2n+k-l]\setminus H\cup S)_o \cup S)
\M(\bar{R}_{n,2n+k-l}([2n+k-l]\setminus H\cup S)_e \cup S).\tag\eai
$$
Indeed, in the special case when all vertices on the symmetry axis $\ell$ are in the same bipartition class (as it is the case for $AR_{2n,2n+k-l}(H,S)$), the factorization theorem of \cite{\FT} states that
$$
\M(G)=2^{w}\M(G^+)\M(G^-)\tag\eaj
$$
where $w$ is half the number of the vertices on $\ell$, and the graphs $G^+$ and $G^-$ are the subgraphs of $G$ obtained above and below $\ell$ by deleting the edges of $G$ incident from above to the 1st, 3rd, 5th, etc. vertex on the symmetry axis (counting from left to right), and deleting the edges incident from below to the 2nd, 4th, 6th, etc. vertex on the symmetry axis. As the set of vertices of $AR_{2n,2n+k-l}(H,S)$ consists of the vertices on $\ell$ with labels in $[2n+k-l]\setminus H\cup S$, (\eai) follows from (\eaj). Using (\eah) twice on the right hand side of (\eai) one obtains~(\eac). \epf

\medskip
\flushpar
{\smc Remark 2.} We presented the case $S=\emptyset$ of Theorem {\taa} in \cite{\FT}.

\mysec{2. The correlation of defects. Statement of the main result}

%
%
%
%
%
%
%
%
%
%

%
%

The correlation of two monomers in a sea of dimers was introduced by Fisher and Stephenson \cite{\FS} as follows. Let $G_{2n}$ be the subgraph of the grid graph $\Z^2$ induced by the vertices contained in the square $[-n,n)\times[n,n)$, and denote by $m_{p,q}$ the monomer consisting of the vertex of $G_{2n}$ of coordinates $(p,q)$. Then the Fisher-Stephenson definition of the correlation $\omega(m_{p,q},m_{p',q'})$ of two oppositely colored (in the chessboard coloring of $G_{2n}$) monomers is
$$
\omega(m_{p,q},m_{p',q'}):=\lim_{n\to\infty}\frac{\M(G_{2n}\setminus\{m_{p,q},m_{p',q'}\})}{\M(G_{2n})}.\tag\ebba
$$
Hartwig \cite{\Hart} proved that\footnote{ The reason for the second form in (\ebbb) is to make the connection with the constant in (\ebad) evident.}
$$
\omega(m_{0,0},m_{d,d-1})\sim\frac{\sqrt{e}}{2^\tfrac56 A^6}\frac{1}{\sqrt{d}}
=\frac{\sqrt{e}}{2^\tfrac{7}{12} A^6}(d\sqrt{2})^{-\tfrac{1}{2}}
,\ \ \ d\to\infty,\tag\ebbb
$$
where 
$$
A=1.28242712...\tag\ebbc
$$ 
is the Glaisher-Kinkelin constant\footnote{ One way to define the Glaisher-Kinkelin constant (see \cite{\Glaish}) is by the limit 
$\lim_{n\to\infty}
\dfrac
 {0!\,1!\,\cdots\,(n-1)!}
 {n^{\frac{n^2}{2}-\frac{1}{12}}\,(2\pi)^{\frac{n}{2}}\,e^{-\frac{3n^2}{4}}}
=
\dfrac
 {e^{\frac{1}{12}}}
 {A}
$.} (eq. (\ebbb), with multiplicative constant given to five decimal places, was conjectured by Fisher and Stephenson in \cite{\FS}).

To make use of the explicit formula given in Theorem {\taa}, it will be convenient for us to define the correlation $\bar\omega$ of defects (monomers and separations) by a variation of the definition of Fisher and Stephenson, namely by enclosing the defects in large Aztec diamonds (the Aztec diamond $AD_n$ is the special case $m=n$ of the Aztec rectangle $AR_{m,n}$). Our definition of correlation applies more generally to any finite set of monomers and separations. However, in view of the set-up of Theorem {\taa}, we will focus on the case when all monomers and separations are along a fixed lattice diagonal.

Throughout this paper we consider the square lattice drawn in the plane so that the lattice lines form angles of degrees $\pm\tfrac{\pi}{4}$ with the horizontal, and we consider a fixed horizontal lattice diagonal $\ell$. For convenient reference to the lattice points on $\ell$, we regard $\ell$ as the number line with unit equal to $\sqrt{2}$ times the lattice spacing, so that the lattice points on $\ell$ are naturally labeled by the integers.

For finite sets of integers $H=\{h_1,\dotsc,h_k\}$ and $S=\{s_1,\dotsc,s_l\}$, $H\cap S=\emptyset$, define the {\it joint correlation $\bar\omega(H;S)$ of having unit holes on $\ell$ at the elements of $H$ and separations on $\ell$ at the elements of $S$} as follows. Our definition is inductive on $|k-l|$:

\medskip
$(i)$ For $k=l$ define 
$$
\bar\omega(h_1,\dotsc,h_k;s_1,\dotsc,s_k)=
\lim_{n\to\infty}\frac
{\M(AD_{2n}(\{h_1,\dotsc,h_k\},\{s_1,\dotsc,s_k\}))}
{\M(AD_{2n})},\tag\ebaa
$$
where the Aztec diamonds on the right hand side are translated so that they are symmetric about $\ell$, and their vertices on $\ell$ have coordinates $-n,-n+1,\dotsc,n-1$

\medskip
$(ii)$ If $k>l$, define

$$
\bar\omega(h_1,\dotsc,h_k;s_1,\dotsc,s_l)=
\frac{1}{\bar\omega\left({\vee \atop \wedge}\right)}
\lim_{d\to\infty}(d\sqrt{2})^{\tfrac{k-l}{2}}\,\bar\omega(a_1,\dotsc,a_k;b_1,\dotsc,b_l,d)
\tag\ebab
$$

\medskip
$(iii)$ If $k<l$, define

$$
\bar\omega(h_1,\dotsc,h_k;s_1,\dotsc,s_l)=
\frac{1}{\bar\omega\left(\circ\right)}
\lim_{d\to\infty}(d\sqrt{2})^{-\tfrac{k-l}{2}}\,\bar\omega(a_1,\dotsc,a_k,d;b_1,\dotsc,b_l),
\tag\ebac
$$
where the values $\bar\omega\left(\circ\right)$ of the correlation of a unit hole and $\bar\omega\left({\vee \atop \wedge}\right)$ of the correlation of a separation are given by
$$
\spreadlines{3\jot}
\align
\bar\omega\left(\circ\right)&:=\frac{e^{\tfrac14}}{2^{\tfrac{7}{24}}A^3}\tag\ebad
\\
\bar\omega\left({\textstyle {\vee \atop \wedge}}\right)&:=\frac{2^{\tfrac{5}{24}}e^{\tfrac14}}{A^3}.
\tag\ebae
\endalign
$$

A {\it cluster of defects} (or defect cluster) is an arbitrary finite union of unit holes and separations on the lattice diagonal $\ell$. Given disjoint clusters of defects $O_1,\dotsc,O_m$ we define their correlation by 
$$
\bar\omega(O_1,\dotsc,O_m):=\bar\omega(H_{O_1\cup\cdots\cup O_m};S_{O_1\cup\cdots\cup O_m}),\tag\ebb
$$
where $H_{O_1\cup\cdots\cup O_m}$ and $S_{O_1\cup\cdots\cup O_m}$ are the set of coordinates of the holes and separations in $O_1\cup\cdots\cup O_m$, respectively.

Note that the correlation of a single defect cluster (the ``self-correlation'') makes sense according to this definition.

In our previous results \cite{\sc} and \cite{\ec} on the asymptotics of the correlation of gaps on the hexagonal lattice, a crucial role was played by the charge $\q(O)$ of a gap $O$, defined to be the difference between the number of white and black monomers in $O$ (in a fixed black and white bipartite coloring of the vertices of the hexagonal lattice). The same turns out to be true for the current situation of the square lattice. 

As indicated by Figure {\faa} and its explanation in the text, each separation defect can be replaced by four trimer gaps. Furthermore, if the vertices on $\ell$ are white, each such trimer contains one white and two black vertices, for an overall charge of $-1$. This and the previous paragraph point to the following extension of the notion of charge to defect clusters:
Given a defect cluster $O$, define its {\it charge} $\q(O)$ to be the number of unit holes in $O$ minus the number of separations in $O$.

The main result of this paper can be stated as follows. 

\proclaim{Theorem \tba}
If $O_1,\dotsc,O_m$ are arbitrary defect clusters on $\ell$, then for large mutual separations between the $O_i$'s, the asymptotics of their correlation is given by
$$
\bar\omega(O_1,\dotsc,O_m)\sim \prod_{i=1}^m \bar\omega(O_i) \prod_{1\leq i<j\leq m} \de(O_i,O_j)^{\frac12 \q(O_i)\q(O_j)},\tag\ebc
$$
where $\de$ is the Euclidean distance.

\endproclaim

We can now explain why we made the indicated choices in $(ii)$--$(iii)$ above when we defined $\bar\omega$ inductively. In $(ii)$, the factor $d\sqrt{2}$ is (asymptotically) the Euclidean distance between the separation at $d$ and the cluster $O$ consisting of the holes at $a_1,\dotsc,a_k$ and the separations at $b_1,\dotsc,b_l$; the exponent $\tfrac{k-l}{2}$ to which it is raised is designed to precisely compensate for the decay of $\bar\omega(a_1,\dotsc,a_k;b_1,\dotsc,b_l,d)$ as the $\frac12\q\left({\vee \atop \wedge}\right)\q(O)=-\tfrac{k-l}{2}$th power of the Euclidean distance between the separation being sent to infinity and the cluster $O$. An analogous remark holds for $(iii)$.

The reason we divide in $(ii)$ and $(iii)$ by the correlation of the extra defect introduced to lower the absolute value of the total charge (for the inductive definition), is to give a chance to (\ebc) to hold (in other words, (\ebab) and (\ebac) are special cases of (\ebc)). What is remarkable is that with the choices (\ebad)--(\ebae) for the correlation of our unit defects, (\ebc) holds in general. 

The value of $\bar\omega(\circ)$ follows from Hartwig's result (\ebbb), if we assume that (\ebc) holds for the special case of two monomers (strictly speaking, we also need the assumption that the $(i)$--$(iii)$-style extension of the Fisher-Stephenson correlation agrees with $\bar\omega$ when restricted to a single monomer and to pairs of monomers; this is a reasonable assumption in view of the fact that, according to the results of Cohn, Elkies and Propp \cite{\CEP}, the dimer statistics is undistorted in the scaling limit at the very center of large Aztec diamonds). 

Finally, it follows from our results in this paper (namely, Corollary {\tbd} and Proposition {\tga}) that 
$$
\bar\omega(\circ)\,\bar\omega\left({\vee \atop \wedge}\right)=\dfrac{\sqrt{e}}{2^{\frac{1}{12}}A^6},\tag\ebd
$$ 
so that the value of $\bar\omega\left({\vee \atop \wedge}\right)$ is determined to be the one given in (\ebae).

\medskip
\flushpar
{\smc Remark 3.} 
The square lattice analog of \cite{\ov, Conjecture\ 1} is that for any {\it monomer clusters}\footnote{ A monomer cluster in an arbitrary finite union of vertices (monomers).} $O_1,\dotsc,O_m$ on $\Z^2$, one has that
$$
\bar\omega(O_1,\dotsc,O_m)\sim\prod_{i=1}^m \bar\omega(O_i) \prod_{1\leq i<j\leq m} \de(O_i,O_j)^{\frac12 \q(O_i)\q(O_j)},\tag\ebe
$$
for large separations between the clusters. 

Theorem {\tba} can be interpreted as an averaged version of (\ebe) for a special class of monomer clusters contained in the union of three consecutive lattice diagonals on $\Z^2$ (namely $\ell$ and the lattice diagonals immediately above and below $\ell$).

Indeed, as seen from Figure {\faa}, each separation defect is equivalent to the superposition of four trimer holes, in the sense that for any graph $G$ with some separations, one of them at $v$, one has
$$
\M(G)=\M(G'\setminus \te_1)+\M(G'\setminus \te_2)+\M(G'\setminus \te_3)+\M(G'\setminus \te_4),\tag\ebff
$$
where $G'$ is the graph obtained from $G$ by ``repairing'' the separation defect at $v$, and $\te_1,\dotsc,\te_4$ are the four trimers centered at $v$ shown in Figure {\faa}(b).

Therefore, the left hand side of (\ebc) is equal to a sum of (at most) $4^l$ terms, each being the correlation of $k$ white monomers and $l$ trimers (each containing one white and two black vertices), where $k$ (resp., $l$) is the number of holes (resp., separations) in the union of the $O_i$'s. Thus, from this point of view, (\ebc) proves an averaged version of conjecture (\ebe) for collinear distributions of holes of {\it arbitrary} charge, {\it even or odd}.

\medskip
\flushpar
{\smc Remark 4.} An interesting feature of (\ebc) is that the multiplicative constant has a simple conceptual interpretation: It is the product of the correlations of the individual defect clusters. On the other hand, the logarithm of the main part on the right hand side of (\ebc) is the sum of the 2D Coulomb electrostatic energies of pairs of electrical charges, so we can view this main part as being given by the superposition principle of 2D electrostatics. For these reasons, and due to the fact that it applies for arbitrary defect clusters on $\ell$, one can regard Theorem {\tba} as a strong superposition principle.

\medskip
\flushpar
{\smc Remark 5.} Our definition $(i)$--$(iii)$ of the correlation is clearly not the only possible one. For instance, rather than defining the correlation of two monomers on $\ell$ in two stages (as this is done through $(i)$--$(iii)$), one could define it in one stage by including a new defect cluster $O$ consisting of two separations, and sending $O$ to infinity (of course, then we need to multiply by the correct power of the Euclidean distance between $O$ and our two monomers, and divide by the correlation of $O$). Theorem {\tba} implies that all such possible variations lead to the same value for the correlation.

\mysec{3. Detailed statement of the main result. Outline of proof}

To phrase (\ebc) more precisely, given a defect cluster $O$ on $\ell$ and an integer $x$, define $O(x)$ to be the translation of $O$ along $\ell$ that takes the leftmost element of $O$ to coordinate $x$ on $\ell$. The main result of this paper can then be phrased in detail as follows.

\proclaim{Theorem \tca} Let $O_1,\dotsc,O_m$ be arbitrary defect clusters on $\ell$, and let $x_1<\cdots<x_m$ be real numbers. Then if $(x_i^{(R)})_R$ is a sequence of integers so that $\lim_{R\to\infty} x_i^{(R)}/R=x_i$, $i=1,\dotsc,m$, we have as $R\to\infty$ that
$$
\spreadlines{3\jot}
\align
\bar\omega\left(O_1(x_1^{(R)}),\dotsc,O_m(x_m^{(R)})\right)&\sim \prod_{i=1}^m \bar\omega(O_i) \prod_{1\leq i<j\leq m} \de\left(O_i(x_i^{(R)}),O_j(x_j^{(R)})\right)^{\frac12 \q(O_i)\q(O_j)}\tag\eca
\\
&
\sim
\prod_{i=1}^m \bar\omega(O_i) \prod_{1\leq i<j\leq m} 
\left(\sqrt{2}(Rx_j-Rx_i)\right)^{\frac12 \q(O_i)\q(O_j)},\tag\ecb
\endalign
$$
where $\de$ is the Euclidean distance between the clusters, defined to be the distance between their leftmost defects.
\endproclaim

{\it Outline of proof.} There are three basic ingredients of the proof: a simple formula stating how the correlation is affected when a monomer and a separation swap places (the exactness lemma, see Lemma {\tbb}), a product formula stating how the correlation changes when a single defect moves one unit (the elementary move lemma, see Lemma {\tbc}), and a set of product formulas giving the exact correlation of defect clusters whose supporting sets have special structure (see Proposition {\tfa}). 

The proof, presented in Section 8, is organized in four steps. Using the elementary move lemma we first reduce to the case when each defect cluster consists of defects occupying consecutive sites of $\Z$ (see Step 1 in Section 8). Next, using exactness, we further reduce to what we call the standard case, i.e. the case when in the union of our clusters all the monomers are to the left of all the separations (this is done in Step 2 of the proof). In Step 3 we prove the standard case in the special case when all defect clusters consist of an even number of defects, using the explicit formulas provided by Proposition {\tfa}. We complete the proof in Step 4 by an induction argument on the rank of the collection of clusters (the rank equals $i$ if the leftmost odd cluster is the $i$th from the right). The induction step amounts to verifying an asymptotic equality which, using the elementary move lemma and the explicit formulas of Proposition {\tfa}, boils down to checking the equality between two explicit products of Gamma functions.

\mysec{4. Exactness}

A crucial role in our results is played by the function $E$ defined by
$$
E(a_1,a_2\dotsc,a_k;b_1,b_2\dotsc,b_l):=\frac
{\prod_{1\leq i<j\leq k}|a_i-a_j|^\frac12\prod_{1\leq i<j\leq l}|b_i-b_j|^\frac12}
{\prod_{i=1}^k\prod_{j=1}^l|a_i-b_j|^\frac12},\tag\ebh
$$
which is the multiplicative version of the Coulomb energy (from two dimensional electrostatics) of the system of electrical charges obtained by including a positive unit charge at the location of each monomer on $\ell$, and a negative unit charge at the location of each separation on $\ell$ (see \cite{\sc,\S\,2} for more on the parallel between random tilings with holes and electrostatics).

The first main ingredient in the proof of Theorem {\tca} is the following.
%

\proclaim{Lemma {\tbb} (Exactness)} For any distinct integers $a_i$ and $b_j$, $i=1,\dotsc,k$, $j=1,\dotsc,l$, the change in correlation caused by changing the hole at $a_1$ into a separation, and the separation at $b_1$ into a hole, is given by
$$
\frac{\bar\omega(a_1,a_2\dotsc,a_k;b_1,b_2\dotsc,b_l)}{\bar\omega(b_1,a_2\dotsc,a_k;a_1,b_2,\dotsc,b_l)}
=
\frac{E(a_1,a_2\dotsc,a_k;b_1,b_2\dotsc,b_l)}{E(b_1,a_2\dotsc,a_k;a_1,b_2,\dotsc,b_l)}.\tag\ebf
$$

\endproclaim


\pf Suppose $k\geq l$. We prove (\ebf) by induction on $k-l$. To check the base case $k=l$, consider the ratio
$$
\frac
{\M(AD_{2n}(a_1,a_2,\dotsc,a_k;b_1,b_2,\dotsc,b_k))}
{\M(AD_{2n}(b_1,a_2,\dotsc,a_k;a_1,b_2,\dotsc,b_k))},
$$
where the Aztec diamonds have been translated so that their vertices on $\ell$ are the integers $-n,-n+1,\dotsc,n-1$. We claim that by applying Theorem {\taa} to the numerator and denominator above, one obtains after simplifications that
$$
\frac
{\M(AD_{2n}(a_1,a_2,\dotsc,a_k;b_1,b_2,\dotsc,b_k))}
{\M(AD_{2n}(b_1,a_2,\dotsc,a_k;a_1,b_2,\dotsc,b_k))}
=
\frac
{
\overset n-1 \to {\underset{i=-n \atop i\neq b_1}\to\prod}\de({\vee \atop \wedge}_{b_1},i)
\,\,\cdot
 \frac
 {\tsize \prod_{j=2}^k \de({\vee \atop \wedge}_{b_1},{\vee \atop \wedge}_{b_j})}
 {\tsize \prod_{i=1}^k \de({\vee \atop \wedge}_{b_1},\boxdot_{a_i})}
}
{
\overset n-1 \to {\underset{i=-n \atop i\neq a_1}\to\prod}\de({\vee \atop \wedge}_{a_1},i)
\,\,\cdot
 \frac
 {\tsize \prod_{j=2}^k \de({\vee \atop \wedge}_{a_1},{\vee \atop \wedge}_{b_j})}
 {\tsize \prod_{i=2}^k \de({\vee \atop \wedge}_{a_1},\boxdot_{a_i})\cdot\de({\vee \atop \wedge}_{a_1},\boxdot_{b_1})}
},\tag\ebj
$$
where, to keep the notation suggestive of the labeled strings involved in the statement of Theorem {\taa} (see Figure {\fac}), $\boxdot_a$ denotes a hole at $a$, and  ${\vee \atop \wedge}_a$ a separation at $a$ (here $\de$ is the distance function between the integers on $\ell$).

To see this, note that the only difference between the labeled strings corresponding to the top and bottom of the left hand side above is that the $O$ and $E$ labels in the cell of coordinate $b_1$ at the numerator have moved to cell $a_1$ at the denominator. Note also that the resulting prefactors of (\fac) are equal, and thus simplify out when taking the ratio on the left hand side of (\ebj). Thus in the ratio on the left hand side of (\ebj) it suffices to consider the factors coming from the $\Delta(\Cal O)\Delta(\Cal E)$ part of formula (\fac), and from these factors only the following do not simplify out: $(i)$ at the numerator, factors contributing to $\Delta(\Cal O)$ that measure the distance from ${\vee \atop \wedge}_{b_1}$ to other $O$-labeled cells, and factors contributing to $\Delta(\Cal E)$ that measure the distance from ${\vee \atop \wedge}_{b_1}$ to other $E$-labeled cells, and $(ii)$ at the denominator, factors contributing to $\Delta(\Cal O)$ that measure the distance from ${\vee \atop \wedge}_{a_1}$ to other $O$-labeled cells, and factors contributing to $\Delta(\Cal E)$ that measure the distance from ${\vee \atop \wedge}_{a_1}$ to other $E$-labeled cells. The factors in $(i)$ represent the distances from $b_1$ to all other integers on $\ell$, with the exceptions that distances to the coordinates of holes are missing, and distances to the coordinates of separations appear twice (once in the $\Delta(\Cal O)$ part and once in the $\Delta(\Cal E)$ part). This explains the numerator on the right hand side of (\ebj). A similar argument explains the denominator, proving (\ebj).

Next, note that the limit as $n\to\infty$ of the ratio of the first products at the numerator and denominator on the right hand side of (\ebj) is 1. Indeed, only $|a_1-b_1|$ factors at the numerator and $|a_1-b_1|$ at the denominator do not simplify out, and the surviving ones are all of the form $n+\alpha$, with $\alpha$ independent of~$n$. Thus the ratio of their products approaches 1 as $n\to\infty$.

Therefore, taking the limit $n\to\infty$ in (\ebj) one obtains
$$
\frac{\bar\omega(a_1,a_2\dotsc,a_k;b_1,b_2\dotsc,b_k)}{\bar\omega(b_1,a_2\dotsc,a_k;a_1,b_2,\dotsc,b_k)}
=
\frac
{\prod_{j=2}^k|b_1-b_j|\prod_{i=2}^k|a_1-a_i|}
{\prod_{i=2}^k|b_1-a_i|\prod_{j=2}^k|a_1-b_j|}.\tag\ebk
$$
One readily checks that this is the same as the right hand side of (\ebf). This completes the proof of the base case $k=l$.

We now check the induction step. For $k>l$, we have by our definition (\ebab) of the correlation $\bar\omega$ that
$$
\frac{\bar\omega(a_1,a_2\dotsc,a_k;b_1,b_2\dotsc,b_l)}{\bar\omega(b_1,a_2\dotsc,a_k;a_1,b_2,\dotsc,b_l)}
=
\lim_{d\to\infty}
\frac{\bar\omega(a_1,a_2\dotsc,a_k;b_1,b_2\dotsc,b_l,d)}{\bar\omega(b_1,a_2\dotsc,a_k;a_1,b_2,\dotsc,b_l,d)}.\tag\ebkk
$$
By the induction hypothesis, the fraction on the right hand side above can be expressed using (\ebh). Using this, one readily checks that taking the limit as $d\to\infty$ on the right hand side of (\ebkk) one obtains the right hand side of (\ebh). This completes the proof by induction of the case $k\geq l$.

The case $k\leq l$ is proved analogously. \epf

\flushpar
{\smc Remark 6.} According to the statement of Theorem {\tca}, we expect (\ebf) to hold {\it in the limit of large separations} between the defects. The above result shows the unexpected fact that (\ebf) holds in fact {\it exactly}, for any {\it finite} separations between the defects. This is why we call the above lemma exactness.

\mysec{5. Elementary move}

The next result gives the change in correlation caused by moving a single defect one unit to the left. Let our system of defects consist of monomers at $a_1,\dotsc,a_k$ and separations at $b_1,\dotsc,b_l$, and set $$
D:=\{a_1,\dotsc,a_k,b_1,\dotsc,b_l\}.
$$
It will be convenient to have a notation for the set of integers strictly in between two distinct integers $x$ and $y$; denote it by $\langle x, y\rangle$. 
The change in correlation under an elementary move turns out to be expressible in terms of the following two ``kernels.'' For any distinct integers $x$ and $y$ define the {\it likes kernel} $L_D(x,y)$ and the {\it unlikes kernel} $U_D(x,y)$ by
$$
L_D(x,y):=
\spreadmatrixlines{3\jot}
\left\{\matrix
\dfrac{\Gamma\left(\frac{|x-y|-1}{2}\right)\Gamma\left(\frac{|x-y|+1}{2}\right)}
{\Gamma^2\left(\frac{|x-y|}{2}\right)},\text{\rm \ \ \ if $|\langle x,y\rangle\setminus D|$ is even}
\\
\dfrac{\Gamma\left(\frac{|x-y|}{2}\right)\Gamma\left(\frac{|x-y|}{2}+1\right)}
{\Gamma^2\left(\frac{|x-y|+1}{2}\right)},\text{\rm \ \ \ if $|\langle x,y\rangle\setminus D|$ is odd}
\endmatrix
\right.
\tag\ebl
$$
and
$$
U_D(x,y):=
\spreadmatrixlines{3\jot}
\left\{\matrix
\dfrac{\Gamma\left(\frac{|x-y|}{2}\right)\Gamma\left(\frac{|x-y|}{2}+1\right)}
{\Gamma^2\left(\frac{|x-y|+1}{2}\right)},\text{\rm \ \ \ if $|\langle x,y\rangle\setminus D|$ is even}
\\
\dfrac{\Gamma\left(\frac{|x-y|-1}{2}\right)\Gamma\left(\frac{|x-y|+1}{2}\right)}
{\Gamma^2\left(\frac{|x-y|}{2}\right)},\text{\rm \ \ \ if $|\langle x,y\rangle\setminus D|$ is odd}
\endmatrix
\right.
\tag\ebm
$$

\proclaim{Lemma {\tbc} (Elementary Move)} Let $a_1,\dotsc,a_k,b_1,\dotsc,b_l$
be distinct integers. 

Then if $a_i-1\notin\{a_1,\dotsc,a_k,b_1,\dotsc,b_l\}$, we have
$$
\frac
{\bar\omega(a_1,\dotsc,a_{i-1},a_i,a_{i+1},\dotsc,a_k;b_1,\dotsc,b_l)}
{\bar\omega(a_1,\dotsc,a_{i-1},a_i-1,a_{i+1},\dotsc,a_k;b_1,\dotsc,b_l)}
=
\frac
{\prod_{j:a_j<a_i}L_D(a_i,a_j)\prod_{j:b_j<a_i}U_D(a_i,b_j)}
{\prod_{j:a_j>a_i}L_D(a_i-1,a_j)\prod_{j:b_j>a_i}U_D(a_i-1,b_j)}
.\tag\ebn
$$

Similarly, if $b_i-1\notin\{a_1,\dotsc,a_k,b_1,\dotsc,b_l\}$, we have
$$
\frac
{\bar\omega(a_1,\dotsc,a_k;b_1,\dotsc,b_{i-1},b_i,b_{i+1},\dotsc,b_l)}
{\bar\omega(a_1,\dotsc,a_k;b_1,\dotsc,b_{i-1},b_i-1,b_{i+1},\dotsc,b_l)}
=
\frac
{\prod_{j:b_j<b_i}L_D(b_i,b_j)\prod_{j:a_j<b_i}U_D(b_i,a_j)}
{\prod_{j:b_j>b_i}L_D(b_i-1,b_j)\prod_{j:a_j>b_i}U_D(b_i-1,a_j)}
.\tag\ebo
$$

\endproclaim

To prove this result, we define an auxiliary correlation $\tilde\omega$ as follows. 

For finite sets of integers $H=\{h_1,\dotsc,h_k\}$ and $S=\{s_1,\dotsc,s_l\}$, with $H\cap S=\emptyset$, define the {\it joint correlation $\tilde\omega$ of having holes on $\ell$ at the elements of $H$ and separations on $\ell$ at the elements of $S$} by
$$
\tilde\omega(H;S)=\tilde\omega(h_1,\dotsc,h_k;s_1,\dotsc,s_l)=
\lim_{n\to\infty}\frac
{\M(AR_{2n,2n+k-l}(\{h_1,\dotsc,h_k\},\{s_1,\dotsc,s_l\}))}
{\M(AR_{2n,2n+k-l}(\{0,1,\dotsc,k-1\},\{k,k+1,\dotsc,k+l-1\}))},\tag\eea
$$
where the Aztec rectangles on the right hand side are translated so that they are symmetric about $\ell$, and their vertices on $\ell$ have coordinates $-n,-n+1,\dotsc,n+k-l-1$. 

First, we prove the statement of Lemma {\tbc} when $\bar\omega$ is replaced by $\tilde\omega$.

\proclaim{Lemma \tbcp} Let $a_1,\dotsc,a_k,b_1,\dotsc,b_l$
be distinct integers. Then if $a_i-1\notin\{a_1,\dotsc,a_k,b_1,\dotsc,b_l\}$, we have
$$
\frac
{\tilde\omega(a_1,\dotsc,a_{i-1},a_i,a_{i+1},\dotsc,a_k;b_1,\dotsc,b_l)}
{\tilde\omega(a_1,\dotsc,a_{i-1},a_i-1,a_{i+1},\dotsc,a_k;b_1,\dotsc,b_l)}
=
\frac
{\prod_{j:a_j<a_i}L_D(a_i,a_j)\prod_{j:b_j<a_i}U_D(a_i,b_j)}
{\prod_{j:a_j>a_i}L_D(a_i-1,a_j)\prod_{j:b_j>a_i}U_D(a_i-1,b_j)}
.\tag\ebnp
$$

Similarly, if $b_i-1\notin\{a_1,\dotsc,a_k,b_1,\dotsc,b_l\}$, we have
$$
\frac
{\tilde\omega(a_1,\dotsc,a_k;b_1,\dotsc,b_{i-1},b_i,b_{i+1},\dotsc,b_l)}
{\tilde\omega(a_1,\dotsc,a_k;b_1,\dotsc,b_{i-1},b_i-1,b_{i+1},\dotsc,b_l)}
=
\frac
{\prod_{j:b_j<b_i}L_D(b_i,b_j)\prod_{j:a_j<b_i}U_D(b_i,a_j)}
{\prod_{j:b_j>b_i}L_D(b_i-1,b_j)\prod_{j:a_j>b_i}U_D(b_i-1,a_j)}
.\tag\ebop
$$

\endproclaim

\pf 
Enclose the holes at $a_1,\dotsc,a_k$ and the separations at $b_1,\dotsc,b_l$ in the Aztec rectangle\footnote{ Our entire argument works just as well with $AR_{2n,2n+k-l}$ instead of $AR_{4n,4n+k-l}$, but then the index of the Pochhammer symbol in (\eboa) and its analogs is $\lfloor\frac{n}{2}\rfloor$ rather than $n$; in order to keep the notation simpler, we stick to  $AR_{4n,4n+k-l}$ as the enclosing region.} \!\!$AR_{4n,4n+k-l}$ translated so that its vertices on the symmetry axis $\ell$ are $-2n,-2n+1,\dotsc,2n+k-l-1$. 
%
To get used to the reasoning we will employ, we work out first the case of just two monomers, at $a_1$ and~$a_2$. 

Consider first the case when there is an even number of integers between $a_1$ and $a_2$. Then if the left monomer moves to the left one unit, the ratio
$$
\frac{\M(AR_{4n,4n+2}(\{a_1,a_2\},\emptyset)}{\M(AR_{4n,4n+2}(\{a_1-1,a_2\},\emptyset)}
$$
can be recorded pictorially by the fraction in front of the equality sign in Figure {\fba}. This fraction is to be interpreted as follows. Each $E$-$O$ labeled string of cells stands for a $\Delta(\Cal O)\Delta(\Cal E)$ product (just as it was the case in Figure {\fac}). More precisely, the locations of the cells are coordinatized from left to right by $-2n,-2n+1,\dotsc,2n+1$ (except for the four shorter strings, whose rightmost cells have coordinate $2n$); then $\Cal O$ and $\Cal E$ are the sets formed by the coordinates of the $O$- and $E$-cells in the string, respectively.

With this interpretation, Theorem {\taa} implies that  $\tfrac{\M(AR_{4n,4n+2}(\{a_1,a_2\},\emptyset)}{\M(AR_{4n,4n+2}(\{a_1-1,a_2\},\emptyset)}$ is equal to the ratio in front of the equality sign in Figure {\fba}.

\topinsert
\centerline{\mypic{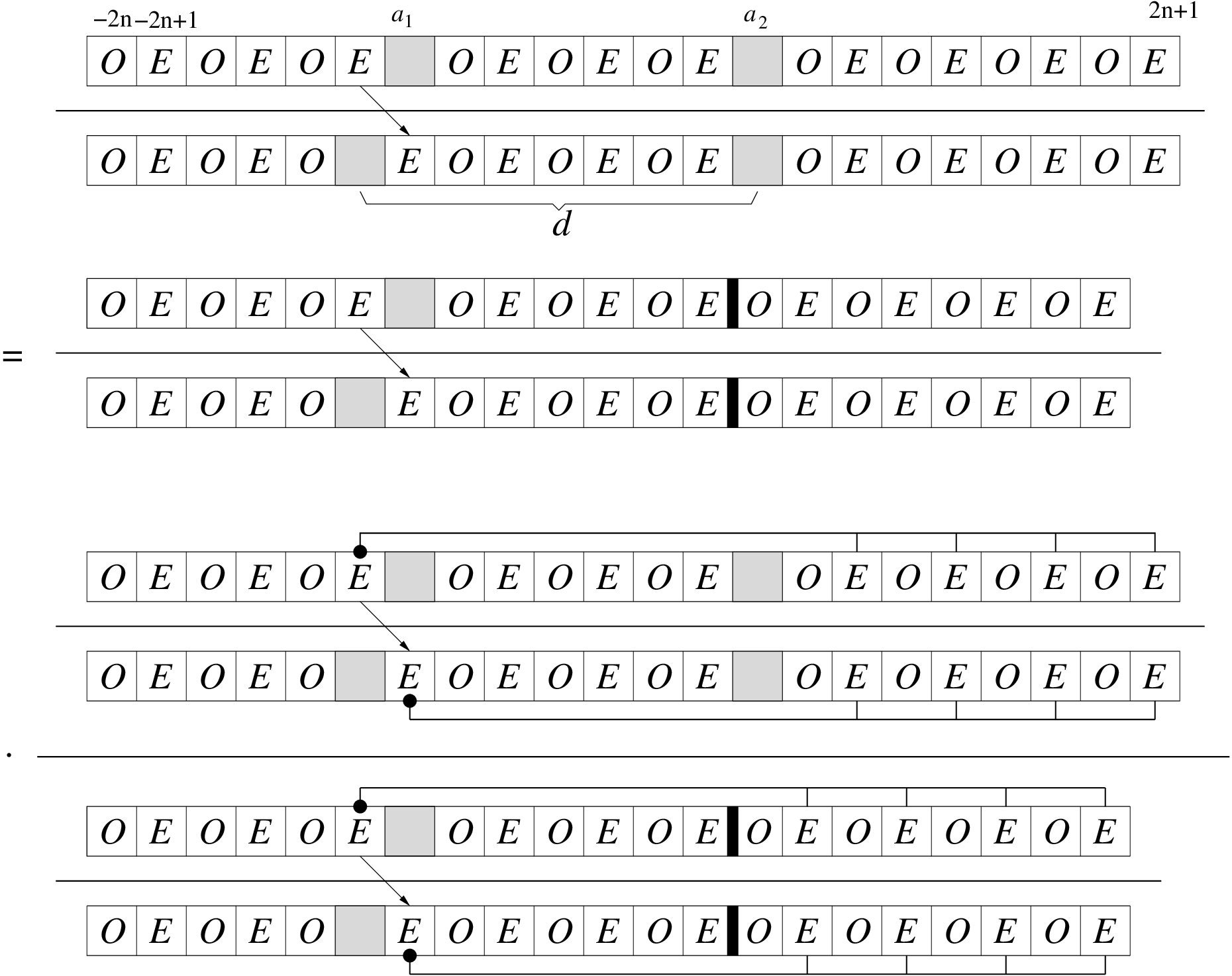}}
\medskip
\centerline{{\smc Figure~{\fba}.} Pictorial interpretation of $\tfrac{\tilde\omega(a_1,a_2)}{\tilde\omega(a_1-1,a_2)}$ when $a_2-a_1-1$ is even; $d=a_2-(a_1-1)$.}
\endinsert

To the right of the equality sign in Figure {\fba} one has the product of two fractions. The first fraction is of the same type as the one in front of the equality sign. In order to explain what the second fraction is, we need to specify the meaning of the decorated strings of cells: Each of them represents the product of the distances from the dotted $E$ to the $E$'s to which it is connected by solid lines.

To see why the equality stated in Figure {\fba} holds, consider how the ratio before the equality sign in the figure changes when the rightmost gaps are eliminated by translating one unit to the left the portions of the strings that are to their right. For both ratios, all factors simplify out except for the ones recording distances from the single ``moving'' $E$ (indicated by arrows in Figure {\fba}) to the other $E$'s. 
After the translation, all leftover factors remain unchanged, except for the the ones representing distances from the moving $E$ to $E$'s in the translated piece. The second ratio after the equality sign precisely corrects this, restoring the equality.


Since the pattern of $O$'s is the same in the numerators and denominators of all fractions in Figure {\fba}, there are no contributions from the $\Delta(\Cal O)$ parts in these fractions.

We claim that the $\Delta(\Cal E)$ part in the first ratio on the right hand side of the equality in Figure {\fba} has the same $n\to\infty$ asymptotics as
$$
\frac{2\cdot4\cdot\dotsc\cdot(2n)}{3\cdot5\cdot\dotsc\cdot(2n+1)}
\cdot
\frac{3\cdot5\cdot\dotsc\cdot(2n+1)}{2\cdot4\cdot\dotsc\cdot(2n)}
=1.\tag{\eboo}
$$
Indeed, the $\Delta(\Cal E)$ part in the first ratio is equal to the left hand side above times a number $c$ of factors of type $\frac{n+\beta}{n+\beta'}$, with $c$, $\beta$, and $\beta'$ independent of $n$.

On the other hand, letting $d=a_2-a_1+1$, the $\Delta(\Cal E)$ part in the second ratio is seen in the same way to have the same $n\to\infty$ asymptotics as
$$
\dfrac
{\frac{(d+2)(d+4)\cdots(d+2n)}{(d+1)(d+3)\cdots(d+2n-1)}}
{\frac{(d+1)(d+3)\cdots(d+2n-1)}{d(d+2)\cdots(d+2n-4)}}
=
\dfrac{\left(\frac{d}{2}\right)_n\left(\frac{d}{2}+1\right)_n}
{\left(\frac{d+1}{2}\right)_n^2},\tag\eboa
$$
where for $k\geq0$ the Pochhammer symbol $(a)_k$ is defined to be
$$
(a)_k:=a(a+1)\cdots(a+k-1).\tag\ebp
$$
Expressing the Pochhammer symbols in terms of Gamma functions by $(a)_k=\Gamma(a+k)/\Gamma(a)$ and using Stirling's formula, one readily sees that 
$$
\lim_{n\to\infty}\frac{\left(\frac{d}{2}\right)_n\left(\frac{d}{2}+1\right)_n}{\left(\frac{d+1}{2}\right)_n^2}
=
\frac{\Gamma^2\left(\frac{d+1}{2}\right)}{\Gamma\left(\frac{d}{2}\right)\Gamma\left(\frac{d}{2}+1\right)},\ \ \ n\to\infty.\tag\ebpa
$$

Putting all the above together, it follows from the equality in Figure {\fba} that 
$$
\lim_{n\to\infty}\frac{\M(AR_{4n,4n+2}(\{a_1,a_2\},\emptyset)}{\M(AR_{4n,4n+2}(\{a_1-1,a_2\},\emptyset)}
=
\frac{\Gamma^2\left(\frac{d+1}{2}\right)}{\Gamma\left(\frac{d}{2}\right)\Gamma\left(\frac{d}{2}+1\right)},\tag\ebpaa
$$
which verifies (\ebn) in this case.

\topinsert
\centerline{\mypic{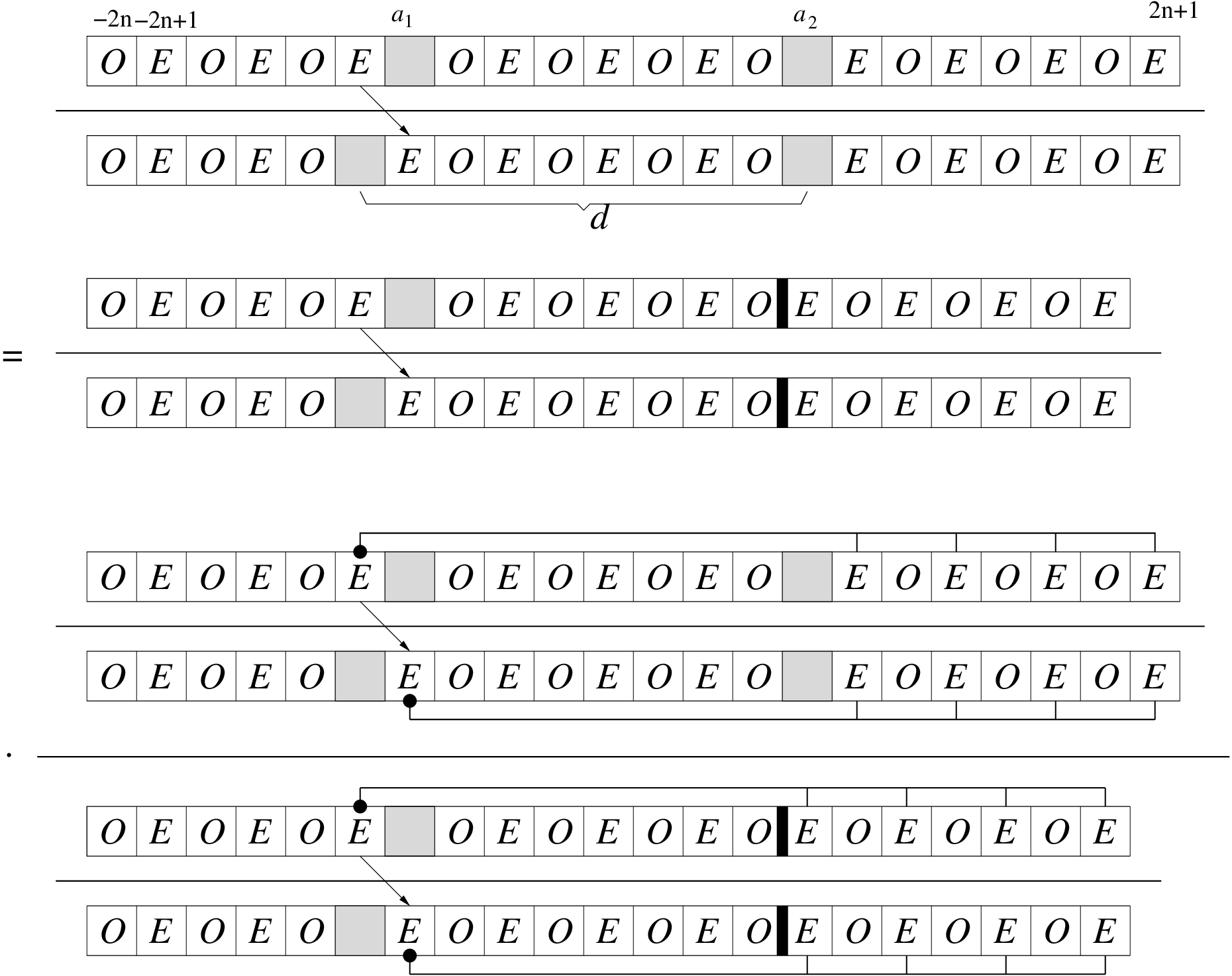}}
\medskip
\centerline{{\smc Figure~{\fbb}.} Interpretation of $\tfrac{\tilde\omega(a_1,a_2)}{\tilde\omega(a_1-1,a_2)}$ when $a_2-a_1-1$ is odd; $d=a_2-(a_1-1)$.}
\endinsert

We now look at how the above calculations change when the number of integers between $a_1$ and $a_2$ is odd. The reason there is a change is that when filling the unaffected cells alternately by $O$'s and $E$'s (as prescribed by Theorem {\taa}), the way the neighborhoods of the two gaps (corresponding to the monomers) look depends on the parity of the number of unaffected cells between the gaps. The case when the latter number is odd is illustrated in Figure {\fbb}. 

Comparing Figure {\fbb} to Figure {\fba} one sees that the patterns of $E$'s they involve are not just similar, but {\it identical}. The only reason there is a difference is because in Figure {\fbb} the value of $d$ is one unit larger that the value of $d$ in Figure {\fba}. Therefore the change in correlation in the current case is obtained by replacing $d$ by $d-1$ in the expression we worked out above in ({\ebpaa}), that is, by
$$
\lim_{n\to\infty}\frac{\M(AR_{4n,4n+2}(\{a_1,a_2\},\emptyset)}{\M(AR_{4n,4n+2}(\{a_1-1,a_2\},\emptyset)}
=
\frac{\Gamma^2\left(\frac{d}{2}\right)}{\Gamma\left(\frac{d-1}{2}\right)\Gamma\left(\frac{d+1}{2}\right)},\tag\ebpaaa
$$
which, as $d=a_2-(a_1-1)$, verifies (\ebn) also for the case when there are an odd number of cells between the gaps.

We are now ready to prove the special case $l=0$ of (\ebn). We proceed by induction on $k$, the number of unit holes on $\ell$. 

For clarity of notation, we present the details of the inductive step from 3 holes to 4 holes. Assume (\ebn) holds for $k=3$, $l=0$, and consider four unit holes at locations $a_1,\dotsc,a_4$. Suppose the hole at $a_2$ is moved one unit to the left. 

Interpreted the same way as Figures {\fba} and {\fbb}, Figure {\fbc} states an equality involving ratios of $\Delta(\Cal O)\Delta(\Cal E)$ products; its proof follows by the same arguments that proved the equalities in Figures {\fba} and~{\fbb}. 
By Theorem {\taa}, this equality is equivalent to
$$
\frac{\M(AR_{4n,4n+4}(\{a_1,a_2,a_3,a_4\};\emptyset)}{\M(AR_{4n,4n+4}(\{a_1,a_2-1,a_3,a_4\};\emptyset)}
=
\frac{\M(AR_{4n,4n+3}(\{a_1,a_2,a_3\};\emptyset)}{\M(AR_{4n,4n+3}(\{a_1,a_2-1,a_3\};\emptyset)}
\cdot F_2,\tag\ebq
$$ 
where $F_2$ is the second fraction on the right hand side of the equality in Figure {\fbc}.

\topinsert
\centerline{\mypic{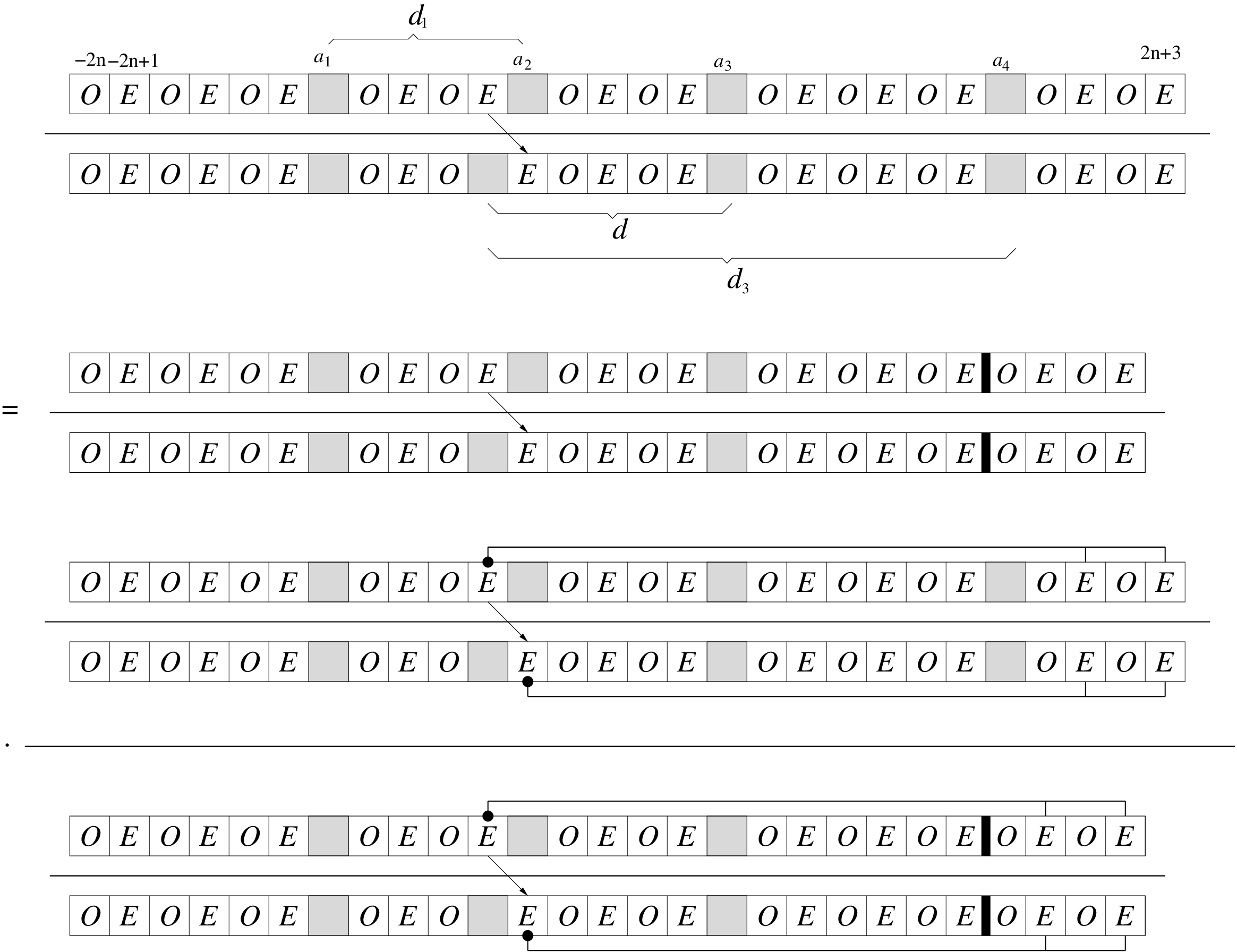}}
\medskip
\centerline{{\smc Figure~{\fbc}.}  Interpretation of $\tfrac{\tilde\omega(a_1,a_2,a_3,a_4)}{\tilde\omega(a_1,a_2-1,a_3,a_4)}$.}
\endinsert

Since the patterns of $O$'s are identical at the numerator and denominator of $F_2$, the $\Delta({\Cal O})$ parts do not contribute to $F_2$. 
For definiteness, assume there is an even number of unaffected cells between $a_2$ and $a_4$. 
By the same arguments we used in the case of two monomers, the $\Delta(\Cal E)$ part of $F_2$ is seen to have the same $n\to\infty$ asymptotics as
$$
\dfrac
{\frac{(d_3+2)(d_3+4)\cdots(d_3+2n)}{(d_3+1)(d_3+3)\cdots(d_3+2n-1)}}
{\frac{(d_3+1)(d_3+3)\cdots(d_3+2n-1)}{d_3(d_3+2)\cdots(d_3+2n-2)}}
=
\dfrac{\left(\frac{d_3}{2}\right)_n\left(\frac{d_3}{2}+1\right)_n}
{\left(\frac{d_3+1}{2}\right)_n^2},\tag\ebr
$$
where $d_3=a_4-(a_2-1)$ (note that this is a calculation not only similar, but {\it identical} to the one in (\eboa); in particular, the fact that there is a gap in between the moving gap and the contracted gap does not affect the calculation at all --- its only effect is that the distance $d_3$ is one unit more than it would be if the intervening gap was not there).

By (\ebr), (\ebpa), and the induction hypothesis, taking the limit as $n\to\infty$ in (\ebq) one obtains
$$
\frac{\tilde\omega(a_1,a_2,a_3,a_4)}{\tilde\omega(a_1,a_2-1,a_3,a_4)}=\frac{L_D(a_2,a_1)}{L_D(a_2-1,a_3)}\frac{1}{L_D(a_2-1,a_4)},\tag\ebrr
$$
thus verifying (\ebn) in this case.

The above calculation assumed, as shown in Figure {\fbc}, that there is an even number of unaffected cells between $a_2$ and $a_4$. The remaining case is handled similarly and leads again to ({\ebrr}), just as~({\ebpaaa}) verified~(\ebn) the same way as ({\ebpaa}) did.

It is clear from the above that the same arguments (with exactly the same calculations) allow going from a system of $k-1$ monomers (out of which one moves one unit to the left) to a system with one extra monomer on the right. To allow for the movement of the rightmost monomer, one parallels the arguments above, but contracting away the {\it leftmost} gap at the induction step. This proves (\ebn) for arbitrary $k$ in the case $l=0$.

\topinsert
\centerline{\mypic{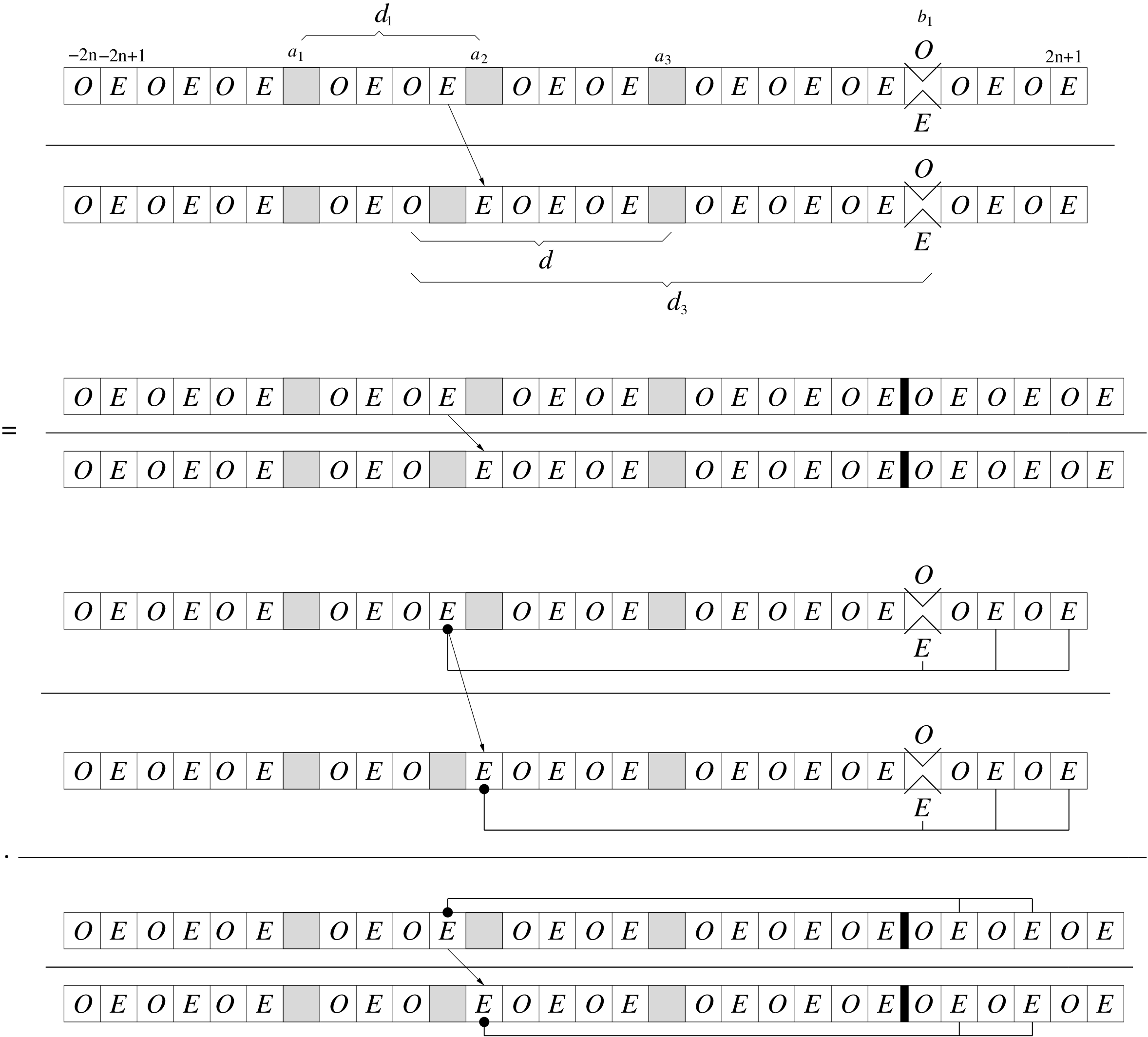}}
\medskip
\centerline{{\smc Figure~{\fbd}.} Interpretation of $\tfrac{\tilde\omega(a_1,a_2,a_3;b_1)}{\tilde\omega(a_1,a_2-1,a_3;b_1)}$.}
\endinsert

To understand the effect on the above reasoning of having also separation defects, we work out the case detailed in Figure {\fbc} with the change that the rightmost defect is a separation instead of a hole.

With our interpretation of $O$-$E$ labeled strings, Figure {\fbd} states an equality whose proof follows by the same arguments we used above. The only difference is that in order to eliminate the overlap at $b_1$ we now {\it pull out} one unit the rightmost piece of the labeled string, instead of pulling it in one unit as in the case of a gap (recall that, as pictured in Figure {\fac}, monomers and separations on $\ell$ correspond to gaps and overlaps, respectively, in the corresponding $O$-$E$ labeled string). 

By Theorem {\taa}, the equality of Figure {\fbd} is equivalent to
$$
\frac{\M(AR_{4n,4n+2}(\{a_1,a_2,a_3\};\{b_1\})}{\M(AR_{4n,4n+2}(\{a_1,a_2-1,a_3\};\{b_1\})}
=
\frac{\M(AR_{4n,4n+3}(\{a_1,a_2,a_3\};\emptyset)}{\M(AR_{4n,4n+3}(\{a_1,a_2-1,a_3\};\emptyset)}
\cdot F'_2,\tag\ebs
$$ 
where $F'_2$ is the second fraction on the right hand side of the equality in Figure {\fbd}.
 
By the same arguments we used in the case of two monomers, $F'_2$ is seen to have the same $n\to\infty$ asymptotics as
$$
\dfrac
{
\frac{(d_3-1)(d_3+1)\cdots(d_3+2n-3)}{(d_3-2)(d_3)\cdots(d_3+2n-4)}
}
{
\frac{(d_3)(d_3+2)\cdots(d_3+2n-2)}{(d_3-1)(d_3+1)\cdots(d_3+2n-3)}
}
=
\dfrac{\left(\frac{d_3-1}{2}\right)^2_n}
{\left(\frac{d_3}{2}-1\right)_n\left(\frac{d_3}{2}\right)_n},\tag\ebt
$$
where now $d_3=b_1-(a_2-1)$ (note again that the only effect of the presence of a gap in between the moving gap and the eliminated overlap is simply that it contributes one extra unit to the distance $d_3$).

Expressing the Pochhammer symbols in terms of Gamma functions as before and using Stirling's formula one readily sees that
$$
\lim_{n\to\infty}\frac{\left(\frac{d-1}{2}\right)_n^2}
{\left(\frac{d}{2}-1\right)_n\left(\frac{d}{2}\right)_n}
=
\frac{\Gamma\left(\frac{d}{2}-1\right)\Gamma\left(\frac{d}{2}\right)}
{\Gamma^2\left(\frac{d-1}{2}\right)},\ \ \ n\to\infty.
$$
Thus, taking the limit $n\to\infty$ in (\ebs) one obtains by the induction hypothesis and (\ebm) that
$$
\frac{\tilde\omega(a_1,a_2,a_3;b_1)}{\tilde\omega(a_1,a_2-1,a_3;b_1)}=\frac{L_D(d_1)}{L_D(d_2)}\frac{1}{U_D(d_3)}
$$
(where $d_1=|a_2-a_1|$, $d_2=|(a_2-1)-a_3|$, and $d_3=|(a_2-1)-b_1|$), thus verifying (\ebn) in this case.

Again, the calculation above assumed (as it is the case in Figure {\fbd}) that there is an even number of unaffected cells between $a_1$ and $b_1$. The remaining case is handled in a similar fashion, and leads again to the above equality. 

It is clear that by these same calculations the truth of (\ebn) for a system $\Cal S$ of holes and separations implies the validity of (\ebn) for the system of defects obtained from $\Cal S$ by including an additional separation to the right of all other defects. 

An analogous argument allows additional separations to be included on the left. Inclusion of additional holes  and separations both on the left and on the right covers the general case of the induction step, proving~(\ebn). Equation (\ebo) is proved in a perfectly similar manner. \epf

{\it Proof of Lemma {\tbc}.} Suppose $k\geq l$. We prove (\ebn) by induction on $k-l$. The base case $k=l$ follows directly from Lemma {\tbcp}, as in this case, by their definitions, $\bar\omega$ and $\tilde\omega$ are the same.

For the induction step, note that by the definition (\ebab) of $\bar\omega$ we have
$$
\frac
{\bar\omega(a_1,\dotsc,a_{i-1},a_i,a_{i+1},\dotsc,a_k;b_1,\dotsc,b_l)}
{\bar\omega(a_1,\dotsc,a_{i-1},a_i-1,a_{i+1},\dotsc,a_k;b_1,\dotsc,b_l)}
=
\lim_{d\to\infty}
\frac
{\bar\omega(a_1,\dotsc,a_{i-1},a_i,a_{i+1},\dotsc,a_k;b_1,\dotsc,b_l,d)}
{\bar\omega(a_1,\dotsc,a_{i-1},a_i-1,a_{i+1},\dotsc,a_k;b_1,\dotsc,b_l,d)}.\tag\ebuone
$$
By the induction hypothesis, the fraction on the right hand side above is given by formula (\ebn). Using this, and the fact that 
$$
\lim_{|x-y|\to\infty}L_D(x,y)=\lim_{|x-y|\to\infty}U_D(x,y)=1,\tag\ebuu
$$ 
(which is readily verified using the defining formulas (\ebl) and (\ebm) and Stirling's formula), one easily sees that the result of taking the limit on the right hand side of (\ebuone) is precisely the expression on the right hand side of (\ebn). This completes the proof of (\ebn) when $k\geq l$. The case $k\leq l$ is handled similarly.

The proof of (\ebo) is analogous. \epf

It will be convenient for the proof of Theorem {\tca} to have specific formulas for the change in correlation, in which one doesn't need to worry about the parity of the number of unaffected cells in various portions of the string. Such formulas are given in Corollary {\tbd} below, which concerns the special case when all runs of unaffected cells between two consecutive defects have even length, and the moving defect moves {\it two} units to the left.

Define modified versions $L(d)$ and $U(d)$ of the kernels (\ebl) and (\ebm) by
$$
L(d):=\frac{\Gamma^2\left(\frac{d-1}{2}\right)\Gamma^2\left(\frac{d+1}{2}\right)}
{\Gamma^4\left(\frac{d}{2}\right)}\tag\ebu
$$
and 
$$
U(d):=\frac{\Gamma^2\left(\frac{d-1}{2}\right)\Gamma^2\left(\frac{d-1}{2}\right)}
{\Gamma\left(\frac{d}{2}-1\right)\Gamma^2\left(\frac{d}{2}\right)\Gamma\left(\frac{d}{2}+1\right)},
\tag\ebv
$$
where $d\geq3$ is an integer.

\proclaim{Corollary \tbd} Let $a_1,\dotsc,a_k$ and $b_1,\dotsc,b_l$ be distinct integers, and denote by $D_<$ the list obtained by sorting the elements of $D=\{a_1,\dotsc,a_k,b_1,\dotsc,b_l\}$ in increasing order. Assume that between any two consecutive elements of $D_<$ there are an even number of integers not in $D$.

Then if $a_i-2\notin\{a_1,\dotsc,a_k,b_1,\dotsc,b_l\}$, we have
$$
\frac
{\bar\omega(a_1,\dotsc,a_{i-1},a_i,a_{i+1},\dotsc,a_k;b_1,\dotsc,b_l)}
{\bar\omega(a_1,\dotsc,a_{i-1},a_i-2,a_{i+1},\dotsc,a_k;b_1,\dotsc,b_l)}
=
\frac
{\prod_{j:a_j<a_i}L(|a_i-a_j|)\prod_{j:b_j<a_i}U(|a_i-b_j|)}
{\prod_{j:a_j>a_i}L(|(a_i-2)-a_j|)\prod_{j:b_j>a_i}U(|(a_i-2)-b_j|)}
.\tag\ebw
$$

Similarly, if $b_i-2\notin\{a_1,\dotsc,a_k,b_1,\dotsc,b_l\}$, we have
$$
\frac
{\bar\omega(a_1,\dotsc,a_k;b_1,\dotsc,b_{i-1},b_i,b_{i+1},\dotsc,b_l)}
{\bar\omega(a_1,\dotsc,a_k;b_1,\dotsc,b_{i-1},b_i-2,b_{i+1},\dotsc,b_l)}
=
\frac
{\prod_{j:b_j<b_i}L(|b_i-b_j|)\prod_{j:a_j<b_i}U(|b_i-a_j|)}
{\prod_{j:b_j>b_i}L(|(b_i-2)-b_j|)\prod_{j:a_j>b_i}U(|(b_i-2)-a_j|)}
.\tag\ebx
$$

\endproclaim

\pf Apply (\ebn) twice to obtain (\ebw), and (\ebo) twice to obtain (\ebx). \epf

\mysec{6. Explicit formulas for the correlation of defects of special support}

In this section we provide explicit simple formulas for the correlation of defect clusters whose supporting set has a special structure. 

The {\it support} $\s(O)$ of a defect cluster $O$ consisting of monomers at $a_1,\dotsc,a_k$ and separations at $b_1,\dotsc,b_l$ is the set $\{a_1,\dotsc,a_k,b_1,\dotsc,b_l\}$. A {\it doublet} is a subset of $\Z$ of the form $\{2s+1,2s+2\}$, with $s\in\Z$. We say that the defect cluster $O$ has {\it nice even support} if $\s(O)$ is the union of doublets. On the other hand, we say that $O$ has {\it nice odd support} if $\s(O)$ is the union of doublets and a singleton
$\{2s+1\}$, with $2s+1$ larger than all the elements in the doublets.

Given a defect cluster $O$ of one of these two special types, we define its {\it canonical rearrangement $\widehat O$} as follows. Let $O$ be a defect cluster of nice even support, and suppose it consists of $k+2i$ monomers and $k$ separations, with $i\geq0$; thus $\s(O)$ is the union of $k+i$ doublets. Then $\widehat O$ is defined to be the cluster consisting of monomers at both positions of the first $i$ doublets (as one scans from left to right), monomers also in the first positions of the remaining $k$ doublets, and separations in the second positions of the latter doublets. In the case when $O$ consists of more separations than monomers, the canonical rearrangement $\widehat O$ is defined the same way, with the only change that at both locations of the initial doublets one places separations rather than monomers.

Let now $O$ be a defect cluster of nice odd support, and assume that it has $k+2i+1$ monomers ($i\geq0$) and $k$ separations. Define the canonical rearrangement $\widehat O$ of this cluster to consist of monomers at both positions of the first $i$ doublets, monomers at the first positions of the remaining $k$ doublets and at the singleton, and separations at the second positions of the last $k$ doublets. In the remaining case when $O$ consists of $k$ monomers and $k+2i+1$ separations with $i\leq0$, make in the definition of $\widehat O$ only the change that at the positions in the first $i$ doublets separations are placed rather than monomers.

%

Given a defect cluster $O$ consisting of monomers at $a_1,\dotsc,a_k$ and separations at $b_1,\dotsc,b_l$, we denote by $E(O)$ the expression $E(a_1,\dotsc,a_k;b_1,\dotsc,b_l)$ (which is defined in (\ebh)).

We are now ready to state the exact formulas we need for the correlation of defect clusters that have nice support.

\proclaim{Proposition \tfa} {\rm (a)}. Suppose that the defect cluster $O$ has nice even support. Then if $O$ consists of $k+2i$ monomers and $k$ separations, $i\geq0$, one has
$$
\bar\omega(O)=\frac{2^\tfrac{i(i-1)}{2}}{\pi^{k+i}}E(\widehat O)E(O),\tag\efa
$$
while if $O$ consists of $k$ monomers and $k+2i$ separations, $i\geq0$, one has
$$
\bar\omega(O)=\frac{2^\tfrac{i(i+1)}{2}}{\pi^{k+i}}E(\widehat O)E(O).\tag\efb
$$

{\rm (b)}. Suppose that the defect cluster $O$ has nice odd support, and that $O$ consists of $k+2i+1$ monomers and $k$ separations. Let the first entries in the doublets of $\s(O)$ be $2s_1+1<\cdots<2s_{k+i}+1$, and let the singleton in $\s(O)$ be $2s+1$. Then if $i\geq0$, we have
$$
\bar\omega(O)=2^\tfrac{i^2}{2}\bar\omega(\circ)\prod_{j=1}^{k+i}\frac{\left(\tfrac12\right)_{s-s_j}^2}{(1)_{s-s_j}(1)_{s-s_j-1}}
E(\widehat O)E(O).\tag\efc
$$
On the other hand, if $i\leq0$, one has
$$
\bar\omega(O)=2^\tfrac{i(i+1)}{2}\bar\omega(\circ)\prod_{j=1}^{k+i}\frac{\left(\tfrac12\right)_{s-s_j}^2}{(1)_{s-s_j}(1)_{s-s_j-1}}
E(\widehat O)E(O).\tag\efd
$$

\endproclaim

The proof of the above Proposition is based on the following exact formula.

\proclaim{Lemma \tfb} For any integers $0\leq s_1<\cdots<s_k\leq n-1$ we have
$$
\spreadlines{4\jot}
\align
&
\frac{\M(AD_{2n}(\{2s_1+1,2s_2+1,\dotsc,2s_k+1\},\{2s_1+2,2s_2+2,\dotsc,2s_k+2\})}
{\M(AD_{2n})}
=
\\
&\ \ \ \ \ \ \ \ \ \ \ \ \ \ \ \ \ \ \ \ \ \ \ \ \ \ \ \ \ 
\prod_{i=1}^k\frac{\left(\tfrac12\right)_{s_i+1}\left(\tfrac12\right)_{n-s_i-1}}{(1)_{s_i}(1)_{n-s_i-1}}
\prod_{1\leq i<j\leq k} \frac{(2s_j-2s_i)^2}{(2s_j-2s_i-1)(2s_j-2s_i+1)}.\tag\efe
\endalign
$$


\endproclaim

\pf Apply Theorem {\taa} to the numerator and denominator on the left hand side of (\efe), but with each element
$i$ of ${\Cal O}$ and ${\Cal E}$ replaced by an indeterminate $a_i$. It is routine to check that after simplifications this yields the expression
$$
\frac{1}{2^k}\prod_{i=1}^k\prod_{j=0 \atop j\neq s_1,\dotsc,s_k}^{n-1}\frac{a_{2s_i+2}-a_{2j+1}}{a_{2s_i+1}-a_{2j+1}}.\tag\eff
$$
Specializing back $a_i=i$ for all indices above one obtains the expression on the right hand side of (\efe). \epf


\proclaim{Corollary \tfc} For any defect cluster $O$ that has nice even support and charge $0$, we have
$$
\bar\omega(O)=\frac{1}{\pi^k}E(\widehat O)E(O),
\tag\efg
$$
where $k$ is the number of monomers (and hence also separations) in $O$.

\endproclaim


\pf First we prove (\efg) in the special case when the defect cluster $O$ consists of monomers at $2s_1+1,\dotsc,2s_k+1$ and separations at $2s_1+2,\dotsc,2s_k+2$ (note that this is equivalent to ${\widehat O}=O$).

Suppose that the left positions in the doublets of $\s(O)$ are at $2s_1+1<\cdots<2s_k+1$. Then, written explicitly, (\efg) for this special case is
$$
\bar\omega(2s_1+1,\dotsc,2s_k+1;2s_1+2,\dotsc,2s_k+2)=\frac{1}{\pi^k}E^2(2s_1+1,\dotsc,2s_k+1;2s_1+2,\dotsc,2s_k+2).
\tag\efh
$$
It is an immediate consequence of Stirling's formula that if $s=\lfloor\tfrac{n}{2}\rfloor+c$, with $c\in\Z$ some constant, then
$$
\lim_{n\to\infty}\frac{\left(\tfrac12\right)_{s+1}\left(\tfrac12\right)_{n-s-1}}{(1)_{s}(1)_{n-s-1}}=\frac{1}{\pi}.
\tag\efi
$$
Together with (\efe) and the definition (\ebaa) of the correlation $\bar\omega$, this implies (\efh). 

Let now $O$ be an arbitrary cluster of charge $0$, having nice even support consisting of $k$ doublets. By~(\efh), we have
$$
\bar\omega(\widehat O)=\frac{1}{\pi^k}E^2(\widehat O).\tag\efj
$$
Dividing side by side equations (\efg) and (\efj), one sees that (\efg) is equivalent to
$$
\frac{\bar\omega(O)}{\bar\omega(\widehat O)}=\frac{E(O)}{E(\widehat O)}.\tag\efk
$$
However, (\efk) follows by exactness (Lemma {\tbb}), and the proof is complete. \epf



To prove prove Proposition {\tfa} we will also need the following auxiliary result.

\proclaim{Lemma \tfe} Defining the correlation of a cluster defect of charge $2$ by including two consecutive separation defects and sending them to infinity yields the same result as we get following our definition $(\ebaa)$-$(\ebac)$ of the correlation $\bar\omega$. More precisely, we have
$$
\align
&
\lim_{d_2\to\infty}\frac
{
(d_2\sqrt{2})^1
\left\{
\lim_{d_1\to\infty}(d_1\sqrt{2})^\frac12
\dfrac{\bar\omega(a_1,\dotsc,a_{k+2};b_1,\dotsc,b_k,d_2,d_1+d_2)}
{\bar\omega\left({\vee \atop \wedge}\right)}
\right\}
}
{
\bar\omega\left({\vee \atop \wedge}\right)
}
=
\\
&\ \ \ \ \ \ \ \ \ \ \ \ \ \ \ \ \ \ \ \ \ \ \ \ \ \ \ \ \ \ \ \ \ \ \ \ \ \ \ \ \ \ \ \ \ \ \ \ \ \ \
\lim_{d\to\infty}\frac{(d\sqrt{2})^2\bar\omega(a_1,\dotsc,a_{k+2};b_1,\dotsc,b_k,d,d+1)}
{\bar\omega\left({\vee \atop \wedge}\!\!{\vee \atop \wedge}\right)}.\tag\efii
\endalign
$$

\endproclaim

\pf We present the details in the case when the defect cluster of charge two consists of two adjacent monomers, say at positions $0$ and $1$ on $\ell$. It will be clear from the proof that the same arguments apply in general.

By Lemma {\tbc} and equations (\ebuu), shifting the position of any of the extra separations by one unit has no effect on the value of the limits in (\efii). Thus it suffices to show that
$$
\align
\lim_{d_2\to\infty}\frac
{
(2d_2\sqrt{2})^1
\left\{
\lim_{d_1\to\infty}(2d_1\sqrt{2})^\frac12
\dfrac{\bar\omega(0,1;2d_2,2d_1+2d_2)}
{\bar\omega\left({\vee \atop \wedge}\right)}
\right\}
}
{
\bar\omega\left({\vee \atop \wedge}\right)
}
=
\lim_{d_2\to\infty}\frac{(2d_2\sqrt{2})^2\bar\omega(0,1;2d_2,2d_2+1)}
{\bar\omega\left({\vee \atop \wedge}\!\!{\vee \atop \wedge}\right)},\tag\efjj
\endalign
$$
which is in turn equivalent to
$$
\lim_{d_2\to\infty}\frac{1}{d_2}
\left\{\lim_{d_1\to\infty}
\sqrt{d_1}\,
\frac{\bar\omega(0,1;2d_2,2d_1+2d_2)}{\bar\omega(0,1;2d_2,2d_2+1)}
\right\}
=
\frac{2^\tfrac34\bar\omega^2\!\!\left({\vee \atop \wedge}\right)}{\bar\omega\left({\vee \atop \wedge}\!\!{\vee \atop \wedge}\right)}.\tag\efkk
$$
It is easy to see, using Corollary {\tbd} and Proposition {\tga}, that the definition (\ebaa)-(\ebac) of the correlation $\bar\omega$ implies
$$
\bar\omega\left(\textstyle{\vee \atop \wedge}\!\!{\vee \atop \wedge}\right)=\frac{2}{\pi}.\tag\efl
$$
Therefore, using the value of $\omega\left({\vee \atop \wedge}\right)$ given by (\ebae), we wee from (\efkk) that what we need to prove is
$$
\lim_{d_2\to\infty}\frac{1}{d_2}
\left\{\lim_{d_1\to\infty}
\sqrt{d_1}\,
\frac{\bar\omega(0,1;2d_2,2d_1+2d_2)}{\bar\omega(0,1;2d_2,2d_2+1)}
\right\}
=
C,
\tag\efm
$$ 
where the constant $C$ is given by (\ege).

We check this as follows. By repeated application of Corollary {\tbd} we obtain that
$$
\frac{\bar\omega(0,1;2d_2,2d_1+2d_2)}{\bar\omega(0,1;2d_2,2d_2+1)}
=
[L(3)]_{d_1-1}[U(2d_2+4)]_{d_1-1}[U(2d_2+5)]_{d_1-1},
\tag\efn
$$
where we used the notation
$$
[f(a)]_k:=f(a)f(a+2)f(a+4)\cdots f(a+2k-2),\tag\efo
$$
for any function $f$ of argument $a$ (recall that $L(d)$ and $U(d)$ are the products of ratios of Gamma functions given by (\ebu) and (\ebv)).

By Proposition {\tga} we obtain that, as $d_1\to\infty$, we have
$$
[L(3)]_{d_1-1}[U(2d_2+4)]_{d_1-1}[U(2d_2+5)]_{d_1-1}\sim
C\sqrt{d_1}\frac{C'_e}{[U(4)]_{d_2-1}}\frac{1}{\sqrt{d_1}}\frac{C'}{[U(3)]_{d_2}}\frac{1}{\sqrt{d_1}},\tag\efp
$$
where the constants $C'$ and $C'_e$ are given by (\egf) and (\egh), respectively.

By (\efn) and (\efp), the inner limit in (\efm) evaluates to
$$
\lim_{d_1\to\infty}
\sqrt{d_1}\,
\frac{\bar\omega(0,1;2d_2,2d_1+2d_2)}{\bar\omega(0,1;2d_2,2d_2+1)}
=
C\frac{C'_e}{[U(4)]_{d_2-1}}\frac{C'}{[U(3)]_{d_2}}.\tag\efq
$$

To prove (\efm) we need the asymptotics of the right hand side above as $d_2\to\infty$. This follows by using Proposition {\tga} again. The convenient thing that happens --- and this happens the same way in the case of a general defect cluster of charge 2 --- is that the asymptotics of the denominators on the right hand side of (\efq) is given by the constant in the corresponding numerators, times $\sqrt{d_2}$. Thus all constants on the right hand side of (\efq) simplify out as $d_2\to\infty$, and the asymptotics is just $Cd_2$. This implies (\efm). \epf

{\it Proof of Proposition {\tfa}.} (a). Suppose $O$ has nice even support, and that $O$ consists of $k+2i$ monomers and $k$ separations, where $i\geq0$. We prove (\efa) by induction on $i$. For $i=0$ the statement follows by Corollary {\tfc}. Let $i>0$, and let the constituents of $O$ be monomers at $a_1,\dotsc,a_{k+2i}$ and separations at $b_1,\dotsc,b_k$. By Lemma {\tfe} we have
$$
\spreadlines{3\jot}
\align
\bar\omega(O)&=\frac{1}{\bar\omega\left({\vee \atop \wedge}\!\!{\vee \atop \wedge}\right)}
\lim_{d\to\infty}(2d\sqrt{2})^{2i}\,\bar\omega(a_1,\dotsc,a_{k+2i};b_1,\dotsc,b_k,2d+1,2d+2)
\\
\ \ \ \ \ \ \ \ \ \ \ \ \ 
&=
\frac{2^i}{\bar\omega\left({\vee \atop \wedge}\!\!{\vee \atop \wedge}\right)}
\lim_{d\to\infty} (2d)^{2i}
\,\bar\omega(a_1,\dotsc,a_{k+2i-1},2d+1;b_1,\dotsc,b_k,a_{k+2i},2d+2)
\\
&\ \ \ 
\times
\frac{E(a_1,\dotsc,a_{k+2i};b_1,\dotsc,b_k,2d+1,2d+2)}{E(a_1,\dotsc,a_{k+2i-1},2d+1;b_1,\dotsc,b_k,a_{k+2i},2d+2)},
\tag\efr
\endalign
$$
where at the second equality we used exactness.

By the induction hypothesis, we have
$$
\spreadlines{3\jot}
\align
&
\bar\omega(a_1,\dotsc,a_{k+2i-1},2d+1;b_1,\dotsc,b_k,a_{k+2i},2d+2)
\\
&\ \ \ \ \ \ \ \ \ \
=
\frac{1}{\pi^{k+i+1}}
E(a_1,\dotsc,a_{k+2i-1},2d+1;b_1,\dotsc,b_k,a_{k+2i},2d+2)
\\
&\ \ \ \ \ \ \ \ \ \ \ \ \ \ \ \ \ \ \ \,
\times
E(2s_1+1,\dotsc,2s_{k+i}+1,2d+1;2s_1+2,\dotsc,2s_{k+i}+2,2d+2),\tag\efs
\endalign
$$
where $\{2s_j+1,2s_j+2\}$, $j=1,\dotsc,k+i$ are the doublets of $\s(O)$.

It readily follows from the definition of $E$ that
$$
\spreadlines{3\jot}
\align
&
\lim_{d\to\infty}E(2s_1+1,\dotsc,2s_{k+i}+1,2d+1;2s_1+2,\dotsc,2s_{k+i}+2,2d+2)
\\
&\ \ \ \ \ \ \ \ \ \ \ \ \ \ \ \ \ \ \ \ \ \
=
E(2s_1+1,\dotsc,2s_{k+i}+1;2s_1+2,\dotsc,2s_{k+i}+2)
\\
&\ \ \ \ \ \ \ \ \ \ \ \ \ \ \ \ \ \ \ \ \ \
=
E(\widehat O),\tag\eft
\endalign
$$
and 
$$
E(a_1,\dotsc,a_{k+2i};b_1,\dotsc,b_k,2d+1,2d+2)
\sim
\frac{1}{(2d)^{2i}}E(a_1,\dotsc,a_{k+2i};b_1,\dotsc,b_k),\ \ \ d\to\infty.\tag\efu
$$
Replacing (\efs), (\eft) and (\efu) in (\efr), and using (\efl), one arrives at (\efa). 

The proof of (\efb) is completely analogous. The reason for the difference in the exponents of 2 on the right hand sides of (\efa) and (\efb) is that, when proving (\efb), at the induction step we divide by $\bar\omega(\circ\circ)$ as opposed to $\bar\omega\left({\vee \atop \wedge}\!\!{\vee \atop \wedge}\right)$, and 
$$
\bar\omega(\circ\circ)=\frac{1}{\pi}\tag\efw
$$ 
(as a simple application of Corollary {\tbd} and Proposition {\tga} shows).

(b). 
%
%
%
%
First we prove (\efc) in the case when $O$ consists of monomers at $2s_1+1,2s_1+2\dotsc,2s_i+1,2s_i+2$ and at $2s_{i+1}+1,\dotsc,2s_{i+k}+1,2s+1$, and separations at  $2s_{i+1}+2,\dotsc,2s_{i+k}+2$ (i.e., when $O=\widehat O$). In view of the definition (\ebae) of the correlation $\bar\omega(O)$, we include an extra separation at location $2d+2$, where $d>s$. By repeated application of Corollary {\tbd} we can gradually bring the extra separation to location $2s+2$, and we obtain that
$$
\spreadlines{3\jot}
\align
&
\frac
{\bar\omega(2s_1+1,2s_1+2,\dotsc,2s_i+1,2s_i+2,2s_{i+1}+1,\dotsc,2s_{i+k}+1,2s+1;2s_{i+1}+2,\dotsc,2s_{i+k}+2,2d+2)}
{\bar\omega(2s_1+1,2s_1+2,\dotsc,2s_i+1,2s_i+2,2s_{i+1}+1,\dotsc,2s_{i+k}+1,2s+1;2s_{i+1}+2,\dotsc,2s_{i+k}+2,2s+2)}
\\
&\ \ \ \ \ \ \ \ \ \ \ \ \ \ \ \ \ \ \ \ \ \ \ \ \ \ \ 
=
[U(3)]_{d-s}
\prod_{j=i+1}^{i+k}\,[U(2s-2s_j+3)]_{d-s}[L(2s-2s_j+2)]_{d-s}
\\
&\ \ \ \ \ \ \ \ \ \ \ \ \ \ \ \ \ \ \ \ \ \ \ \ \ \ \ \ \ \ \ \ \ \ \ \ \ \ \ \ \ \,
\times
\prod_{j=1}^i\,[U(2s-2s_j+3)]_{d-s}[U(2s-2s_j+2)]_{d-s}
\tag\efmm
\endalign
$$
(recall that $[f(a)]_k$ is defined by (\efo)).

The $d\to\infty$ asymptotics of the right hand side of (\efmm) is 
$$
\spreadlines{3\jot}
\align
&
[U(3)]_{d-s}\prod_{j=i+1}^{i+k}\,[U(2s-2s_j+3)]_{d-s}[L(2s-2s_j+2)]_{d-s}
\prod_{j=1}^i\,[U(2s-2s_j+3)]_{d-s}[U(2s-2s_j+2)]_{d-s}
\\
&\ \ \ \ \ \ \ \ \ \ \ \ \ 
=
[U(3)]_{d-s}
\prod_{j=i+1}^{i+k}\,\frac{[U(3)]_{d-s_j}}{[U(3)]_{s-s_j}}\frac{[L(4)]_{d-s_j-1}}{[L(4)]_{s-s_j-1}}
\prod_{j=1}^{i}\,\frac{[U(3)]_{d-s_j}}{[U(3)]_{s-s_j}}\frac{[U(4)]_{d-s_j-1}}{[U(4)]_{s-s_j-1}}
\\
&\ \ \ \ \ \ \ \ \ \ \ \ \ 
\sim
\frac{C'}{\sqrt{d}}
\prod_{j=i+1}^{i+k} \frac{\dfrac{C'}{\sqrt{d}}}{[U(3)]_{s-s_j}} \frac{C_e\sqrt{d}}{[L(4)]_{s-s_j-1}}
\prod_{j=1}^{i} \frac{\dfrac{C'}{\sqrt{d}}}{[U(3)]_{s-s_j}} \frac{\dfrac{C'_e}{\sqrt{d}}}{[U(4)]_{s-s_j-1}}
\\
&\ \ \ \ \ \ \ \ \ \ \ \ \ 
=
\frac{2^{k+i-1}C}{\pi^{k+i}d^{i+\frac12}}\prod_{j=i+1}^{k+i} \frac{1}{[U(3)]_{s-s_j}[L(4)]_{s-s_j-1}}
\prod_{j=1}^{i} \frac{1}{[U(3)]_{s-s_j}[U(4)]_{s-s_j-1}},\tag\efoo
\endalign
$$
where we used Proposition {\tga} at the second step, and the facts that $C'=\dfrac{C}{2}$, $C'_e=C_e$, and $CC_e=\dfrac{4}{\pi}$ ($C$, $C'$, $C_e$ and $C'_e$ are the multiplicative constants in the statement of Proposition {\tga}). By part (a), it follows from~(\efmm) and (\efoo) that 
$$
\spreadlines{3\jot}
\align
&
\bar\omega(\widehat O)=\bar\omega(2s_1+1,2s_1+2,\dotsc,2s_i+1,2s_i+2,2s_{i+1}+1,\dotsc,2s_{i+k}+1,2s+1;2s_{i+1}+2,\dotsc,2s_{i+k}+2)
\\
&
=
\frac{1}{\bar\omega\left({\vee \atop \wedge}\right)}
\lim_{d\to\infty}(2d\sqrt{2})^{i+\frac12}
\\
& 
\times
\bar\omega(2s_1+1,2s_1+2,\dotsc,2s_i+1,2s_i+2,2s_{i+1}+1,\dotsc,2s_{i+k}+1,2s+1;2s_{i+1}+2,\dotsc,2s_{i+k}+2,2d+2)
\\
&
= 
\frac{2^{k+\frac{i(i+1)}{2}-\frac14}C}{\pi^{2k+2i}\bar\omega\left({\vee \atop \wedge}\right)}
E^2(2s_1+1,\dotsc,2s_{k+i}+1,2s+1;2s_1+2,\dotsc,2s_k+2,2s+2)
\\
&\ \ 
\times
\prod_{j=i+1}^{k+i} \frac{1}{[U(3)]_{s-s_j}[L(4)]_{s-s_j-1}}
\prod_{j=1}^{i} \frac{1}{[U(3)]_{s-s_j}[U(4)]_{s-s_j-1}}
.\tag\efpp
\endalign
$$
Using the value of $\bar\omega\left({\vee \atop \wedge}\right)$ given by (\ebae), the value of the constant $C$ given by (\ege), equations (\egk) and (\egl), and the definition (\ebu) of $L(d)$, it is routine to check that the right hand side of (\efpp) agrees with the right hand side of (\efc) specialized to the case $O=\widehat O$.

The case of an arbitrary cluster $O$ of nice odd support, containing $k+2i$  follows now easily, using the same idea as in the proof of Corollary {\tfc}. Indeed, suppose $\s(O)$ consists of the $k$ doublets $\{2s_1+1,2s_1+2\},\dotsc,\{2s_k+1,2s_k+2\}$ and the singleton $\{2s+1\}$. By the case we have just proved, we have
$$
\bar\omega(\widehat O)=\bar\omega\left({\textstyle{\vee \atop \wedge}}\right)\prod_{j=1}^{k+i} \frac{\left(\frac12\right)_{s-s_j}^2}{(1)_{s-s_j}(1)_{s-s_j-1}}
\,E^2(\widehat O).
\tag\efqq
$$
By dividing side by side equations (\efc) and (\efqq), one sees that (\efc) is equivalent to
$$
\frac{\bar\omega(O)}{\bar\omega(\widehat O)}=\frac{E(O)}{E(\widehat O)},
$$
which follows by exactness (Lemma {\tbb}). This completes the proof of (\efc). The proof of (\efd) is completely analogous. \epf

\medskip
\flushpar
{\smc Remark 7.} Recall that each separation defect is equivalent to summing over four trimer defects (see Figure~1.1, its explanation in Section 1, and Remark 3), and thus the correlation of a defect cluster consisting of $k$ monomers and $l$ separations is equal to the sum of $4^l$ correlations of $k$ monomers and $l$ trimers (some of these correlations may be 0, if there are consecutive separations).

It is interesting to note that, while the sum of these $4^l$ terms has a nice product expression (by the results in this section and~Remark 8), the individual terms are not so simple in general. Indeed, for instance it turns out that we have$$
\bar\omega(0,1;5,6)=\frac{1}{\pi^2},
$$
while the sixteen correlations of two monomers and two trimers, whose sum $\bar\omega(0,1;5,6)$ is equal to by Remark 3, have the following factorizations:
$$
\align
&
\frac{1755 \pi^3  - 9162 \pi^2  + 22512 \pi - 19136}{16200\pi^4},
\frac{1755 \pi^3  - 9162 \pi^2  + 22512 \pi - 19136}{16200\pi^4},
\frac{26(3\pi-4)(15\pi-26)}{2025\pi^4},
\\
&
-\frac{(3 \pi - 13) (9 \pi^2  - 24 \pi + 32)}{405\pi^4},
\frac{(45 \pi^2  - 138 \pi + 208) (225 \pi^2  - 1380 \pi + 3016)}{324000 \pi^4},
\\
&
-\frac{15525 \pi^3  - 120240 \pi^2  + 230664 \pi - 126464}{162000 \pi^4},
\frac{3375 \pi^3  - 20880 \pi^2  + 64272 \pi - 62192}{40500\pi^4},
\\
&
-\frac{(45 \pi^2  - 168 \pi + 368) (45 \pi^2  - 228 \pi + 208)}{64800\pi^4},
-\frac{15525 \pi^3  - 120240 \pi^2  + 230664 \pi - 126464}{162000 \pi^4},
\\
&
\frac{(45 \pi^2  - 138 \pi + 208) (225 \pi^2  - 1380 \pi + 3016)}{324000 \pi^4},
\frac{3375 \pi^3  - 20880 \pi^2  + 64272 \pi - 62192}{40500\pi^4},
\\
&
-\frac{(45 \pi^2  - 168 \pi + 368) (45 \pi^2  - 228 \pi + 208)}{64800\pi^4},
-\frac{(45 \pi^2  - 228 \pi + 208) (225 \pi^2  - 1020 \pi + 2392)}{324000\pi^4},
\\
&
-\frac{(45 \pi^2  - 228 \pi + 208) (225 \pi^2  - 1020 \pi + 2392)}{324000\pi^4},
-\frac{(3 \pi - 13) (225 \pi^2  - 780 \pi + 1352)}{10125\pi^4},
\\
&
\frac{(45 \pi^2  - 228 \pi + 208)^2}{32400\pi^4}.
\endalign
$$
The fact that, even though individually such terms do not have simple expressions, their sum is given by a simple product formula, is a key ingredient for the success of the method of proof presented in this paper.

\medskip
\flushpar
{\smc Remark 8.} We note that any defect cluster can be changed into a defect cluster of nice support
by applying successive elementary moves to its constituent monomers and/or separations. If the defect cluster $O$ can be transformed into one with nice support by a small number of elementary moves, the elementary move lemma (Lemma {\tbc}) together with Proposition {\tfa} yield a simple formula for $\bar\omega(O)$. In general, however, constituents may need to be moved long distances in order to achieve nice support (e.g., if the distances between all pairs of consecutive defects is large), and then the resulting formula is of more limited use.

%
%
%
%
%
%

\mysec{7. The asymptotics of products of consecutive $L(i)$'s and $U(i)$'s}

Repeated application of Corollary {\tbd} gives rise to products of $L(i)$'s and $U(i)$'s in which the arguments increase by 2 units from one factor to the next. The asymptotics of such products is given by the following result.

\proclaim{Proposition \tga} As $d\to\infty$ the following asymptotic formulas hold:
$$
\align
\prod_{i=1}^d L(2i+1)&\sim C\sqrt{d}\tag\ega\\
\prod_{i=1}^d U(2i+1)&\sim \frac{C'}{\sqrt{d}}\tag\egb\\
\prod_{i=2}^d L(2i)&\sim C_e\sqrt{d}\tag\egc\\
\prod_{i=2}^d U(2i)&\sim \frac{C'_e}{\sqrt{d}},\tag\egd\\
\endalign
$$
where the multiplicative constants on the right hand sides are given by
$$
\align
C&=\frac{2^\frac16 \sqrt{e}\,\pi}{A^6}\tag\ege\\
C'&=\frac{C}{2}\tag\egf\\
C_e&=\frac{4}{\pi C}\tag\egg\\
C'_e&=C_e\tag\egh\\
\endalign
$$
$($recall that $A$ is the Glaisher-Kinkelin constant $(\ebbc)$$)$.

\endproclaim

\pf By the definition (\ebu) of $L$ we have that the left hand side of (\ega) is
$$
\prod_{i=1}^d L(2i+1)=\prod_{i=1}^d\left(\frac{\Gamma(i)\Gamma(i+1)}{\Gamma^2\left(i+\tfrac12\right)}\right)^2.
$$
Using the recurrence $\Gamma(x+1)=x\Gamma(x)$ and that $\Gamma\left(\tfrac12\right)=\sqrt{\pi}$, one can rewrite
the factor in the parentheses above as
$$
\frac{\Gamma(i)\Gamma(i+1)}{\Gamma^2\left(i+\tfrac12\right)}
=
\frac{2^{4i-1}}{\pi}\frac{i!^2(i-1)!^2}{(2i)!(2i-1)!}.
$$
Therefore, we have
$$
\prod_{i=1}^d L(2i+1)=\left(\prod_{i=1}^{d} \frac{2^{4i-1}}{\pi}\frac{i!^2(i-1)!^2}{(2i)!(2i-1)!}\right)^2
=\left(\frac{2^{2d(d+1)-d}}{\pi^d}\frac{(0!\,1!\cdots(d-1)!)^2(1!\,2!\cdots(d)!)^2}{1!\,2!\cdots(2d)!}\right)^2.
\tag\egi
$$
The asymptotics of the product of successive factorials is given by Glaisher's formula (see \cite{\Glaish}):
$$
0!\,1!\cdots(n-1)!\sim \frac{e^{\frac{1}{12}}}{A}\,
n^{\tfrac{n^2}{2}-\tfrac{1}{12}}\,(2\pi)^{\tfrac{n}{2}}\,e^{-\tfrac{3n^2}{4}},\ \ \ n\to\infty,\tag\egj
$$
where $A$ is the Glaisher-Kinkelin constant (\ebbc).

The asymptotics of the first product of factorials at the numerator of (\egi) follows directly by (\egj). The asymptotics of the second product of factorials at the numerator follows by viewing it as the product $(0!\,1!\cdots(d-1)!)\cdot(d!)$, and using (\egi) and Stirling's formula; the asymptotics of the denominator is determined similarly. Putting all this together, it follows that the asymptotics of the right hand side of (\egi) is given by the expression on the right hand side of (\ega). This proves (\ega). 

It follows from the definition of $L(d)$ and $U(d)$ that
$$
\prod_{i=1}^d U(2i+1)=\frac{\Gamma\left(\frac32\right)\Gamma\left(d+\frac12\right)}{\Gamma\left(\frac12\right)\Gamma\left(d+\frac32\right)}\prod_{i=1}^d L(2i+1)=\frac{1}{2d+1}\prod_{i=1}^d L(2i+1).\tag\egk
$$
Together with (\ega), this implies (\egb).

A calculation similar to the one that proved (\ega) proves (\egc). Since
$$
\prod_{i=2}^d U(2i)=\frac{\Gamma(2)\Gamma(d)}{\Gamma(1)\Gamma(d+1)}\prod_{i=2}^d L(2i)=\frac{1}{d}\prod_{i=2}^d L(2i),\tag\egl
$$
one sees that (\egd) follows from (\egc). \epf

\mysec{8. The proof of Theorem {\tca}}

Given a defect cluster $O$, we define its {\it compression $\com(O)$} to be the defect cluster obtained from $O$ as follows. 
For $a\in\Z$, let $\ho(a)$ and $\se(a)$ denote a unit hole at $a$ and a separation at $a$, respectively. Then if 
$$
O=\De_1(a_1)\cup\De_2(a_2)\cup\cdots\cup\De_t(a_t),
$$
with $a_1<\cdots<a_t$ and $\De_i\in\{\ho,\se\}$, $i=1,\dotsc,t$, we define the compression $\co(O)$ of $O$ by
$$
\co(O):=\De_1(a_1)\cup\De_2(a_1+1)\cup\cdots\cup\De_t(a_1+t-1)
$$
(in other words, $\co(O)$ is obtained from $O$ by translating all its constituent holes and separations so as to form a contiguous run starting at $a_1$).

A defect cluster $O$ is said to be {\it compressed} if $\co(O)=O$. 

{\it Proof of Theorem {\tca}.} Since when we coordinatized the lattice points on $\ell$ by the integers we used a unit of $\sqrt{2}$ (in accordance with the lattice spacing being 1), we have that
$$
\de\left(O_i(x_i^{(R)}),O_i(x_i^{(R)})\right)=\sqrt{2}\left(x_j^{(R)}-x_i^{(R)}\right).
$$
Since the hypothesis implies that $x_j^{(R)}-x_i^{(R)}\sim Rx_j-Rx_i$, as $R\to\infty$, this shows that (\eca) and (\ecb) are equivalent.

To prove (\ecb) we proceed in the following steps.

{\smc Step 1.} {\it Reduction to compressed clusters.} Consider the cluster $O_i(x_i^{(R)})$. Bring all defects in it next to the leftmost one, using Lemma {\tbc}. By (\ebuu), at each elementary move, every factor on the right hand sides of (\ebn) and (\ebo) coming from a defect outside $O_i$ approaches 1 as $R\to\infty$. Therefore, to find the change in the correlation $\bar\omega(O_1,\dotsc,O_m)$ due to the compressing of $O_i$, it is enough to keep track of the contributing factors coming from defects within $O_i$. We claim that the cumulative effect of the latter is multiplication by $\dfrac{\bar\omega(O_i)}{\bar\omega(\co(O_i))}$.

Indeed, according to Lemma {\tbc}, precisely the same factors arise if we consider just the defect cluster $O_i$ by itself, and bring all defects in it next to the leftmost one by elementary moves (so that $O_i$ is transformed into $\co(O_i)$). Thus we proved that
$$
\frac
{\bar\omega\left(O_1(x_1^{(R)}),\dotsc,O_m(x_m^{(R)})\right)}
{\bar\omega\left(O_1(x_1^{(R)}),\dotsc,O_{i-1}(x_{i-1}^{(R)}),\co(O_i(x_i^{(R)})),O_{i+1}(x_{i+1}^{(R)}),\dotsc,O_m(x_m^{(R)})\right)}
=
\frac{\bar\omega(O_i)}{\bar\omega(\co(O_i))}.\tag\eha
$$
Applying (\eha) for $i=1,\dotsc,m$ we obtain
$$
\frac
{\bar\omega\left(O_1(x_1^{(R)}),\dotsc,O_m(x_m^{(R)})\right)}
{\bar\omega\left(\co(O_1(x_1^{(R)})),\dotsc,\co(O_m(x_m^{(R)}))\right)}
=
\frac
{\bar\omega(O_1)\cdots\bar\omega(O_m)}
{\bar\omega(\co(O_1))\cdots\bar\omega(\co(O_m))}.\tag\ehb
$$
Rewriting (\ehb) as
$$
\frac
{\bar\omega\left(O_1(x_1^{(R)}),\dotsc,O_m(x_m^{(R)})\right)}
{\bar\omega(O_1)\cdots\bar\omega(O_m)}
=
\frac
{\bar\omega\left(\co(O_1(x_1^{(R)})),\dotsc,\co(O_m(x_m^{(R)}))\right)}
{\bar\omega(\co(O_1))\cdots\bar\omega(\co(O_m))},\tag\ehc
$$
one sees that it suffices to prove (\ecb) for compressed defect clusters.

\medskip
{\smc Step 2.} {\it Reduction of the compressed case to the standard case.} Let $O_1,\dotsc,O_m$ be a collection of compressed defect clusters, and assume that $O_i$ lies entirely to the left of $O_j$, for any $1\leq i<j\leq m$. We define the {\it standardization} $\st(O_1),\dotsc,\st(O_m)$ of this collection as follows. Suppose that $O_1\cup\cdots\cup O_m$ contains a total of $k$ monomers and $l$ separations. Let $\st(O)$ be the defect cluster obtained by 
placing monomers in the first $k$ positions (from left to right) of $\s(O_1)\cup\cdots\cup\s(O_m)$, and separations in the remaining $l$ positions. Then $\st(O_i)$ is defined to be the restriction of $\st(O)$ to $\s(O_i)$.

By exactness (Lemma {\tbb}), we have
$$
\frac
{\bar\omega\left(O_1(x_1^{(R)}),\dotsc,O_m(x_m^{(R)})\right)}
{\bar\omega\left(\st(O_1(x_1^{(R)})),\dotsc,\st(O_m(x_m^{(R)}))\right)}
=
\frac
{E\left(O_1(x_1^{(R)}),\dotsc,O_m(x_m^{(R)})\right)}
{E\left(\st(O_1(x_1^{(R)})),\dotsc,\st(O_m(x_m^{(R)}))\right)}.\tag\ehd
$$
It readily follows from the definition (\ebh) of $E$ that, as $R\to\infty$, one has
$$
\frac
{E\left(O_1(x_1^{(R)}),\dotsc,O_m(x_m^{(R)})\right)}
{E\left(\st(O_1(x_1^{(R)})),\dotsc,\st(O_m(x_m^{(R)}))\right)}
\sim
\prod_{i=1}^m\frac{E(O_i)}{E(\st(O_i))}
\prod_{1\leq i<j\leq m}\frac{(Rx_j-Rx_i)^{\frac12 \q(O_i)\q(O_j)}}{(Rx_j-Rx_i)^{\frac12 \q(\st(O_i))\q(\st(O_j))}},
\tag\ehe
$$
where we have also used that $x_i^{(R)}\sim Rx_i$ as $R\to\infty$, by the hypothesis of Theorem {\tca}. By using exactness again, (\ehd) and (\ehe) imply
$$
\frac
{\bar\omega\left(O_1(x_1^{(R)}),\dotsc,O_m(x_m^{(R)})\right)}
{\bar\omega\left(\st(O_1(x_1^{(R)})),\dotsc,\st(O_m(x_m^{(R)}))\right)}
\sim
\prod_{i=1}^m\frac{\bar\omega(O_i)}{\bar\omega(\st(O_i))}
\prod_{1\leq i<j\leq m}\frac{(Rx_j-Rx_i)^{\frac12 \q(O_i)\q(O_j)}}{(Rx_j-Rx_i)^{\frac12 \q(\st(O_i))\q(\st(O_j))}}.
\tag\ehf
$$
Therefore, if (\ecb) holds for the collection of defect clusters $\{\st(O_1),\dotsc,\st(O_m)\}$, it follows from (\ehf) that~(\ecb) also holds for the collection $\{O_1,\dotsc,O_m\}$.

We call a defect cluster even or odd according to the parity of the number of defects it consists of.

\medskip
{\smc Step 3.} {\it Proof of the standard case when all clusters are even.} Let $O_1,\dotsc,O_m$ be a standardized collection of compressed even clusters. Then for each $1\leq i\leq m$, $\s(O_i(x_i^{(R)}))$ consists of an even number of consecutive integers. By repeated use of Lemma {\tbc}, we may successively shift each defect cluster one unit, if necessary, and position them so that they all have nice even support (this is desirable, as then we can apply the exact formulas of Proposition {\tfa}). Moreover, by (\ebuu), this will make no change in the $R\to\infty$ asymptotics of their correlation. Therefore we may assume that, for all $i$, $O_i(x_i^{(R)})$ --- and therefore $O_1(x_1^{(R)})\cup\cdots\cup O_m(x_m^{(R)})$ as well --- have nice even support. Then the exact formulas (\efa)-(\efb) apply to all correlations involved in (\ecb), and proving (\ecb) in this case reduces to checking the agreement of the two resulting sides.

Since by assumption our collection of defect clusters is standardized, it follows that there is and index $u$, $1\leq u\leq m$, so that $O_1,\dotsc,O_{u-1}$ consist entirely of monomers,  $O_{u+1},\dotsc,O_{m}$ consist entirely of separations, and $O_u$ consists of a run of say $p_u$ monomers followed by a run of say $n_u$ separations, where $p_u,n_u\geq0$. Let $p_i$ be the number of monomers in $O_i$, for $i=1,\dotsc,u-1$, and let $n_j$ be the number of separations in~$O_j$, for $j=u+1,\dotsc,m$. Set $O^{(R)}:=O_1(x_1^{(R)})\cup\cdots\cup O_1(x_1^{(R)})$, $P:=\sum_{i=1}^u p_i$, and $N:=\sum_{j=u}^m n_j$. To apply Proposition {\tfa}(a), one needs to know whether $P\geq N$ or $P\leq N$, and whether $p_u\geq n_u$ or $p_u\leq n_u$. Suppose $P\geq N$, and $p_u\geq n_u$. Then, by Proposition {\tfa}(a), (\ecb) amounts to
$$
\spreadlines{4\jot}
\align
\frac{2^\frac{\frac{P-N}{2}\left(\frac{P-N}{2}-1\right)}{2}}{\pi^{\frac{P+N}{2}}}
E\left(\widehat{O^{(R)}}\right)E\left({O^{(R)}}\right)
&\sim
\frac{2^\frac{\frac{p_u-n_u}{2}\left(\frac{p_u-n_u}{2}-1\right)}{2}}{\pi^\frac{p_u+n_u}{2}}
E(\widehat{O}_u)E(O_u)
\\
&
\times
\prod_{i=1}^{u-1}\frac{2^\frac{\frac{p_i}{2}\left(\frac{p_i}{2}-1\right)}{2}}{\pi^\frac{p_i}{2}}
E(\widehat{O}_i)E(O_i)
\\
&
\times
\prod_{i=u+1}^{m}\frac{2^\frac{\frac{n_i}{2}\left(\frac{n_i}{2}+1\right)}{2}}{\pi^\frac{n_i}{2}}
E(\widehat{O}_i)E(O_i)
\\
&
\times
\prod_{1\leq i<j\leq m} \left(\sqrt{2}(Rx_j-Rx_i)\right)^{\frac12 \q(O_i)\q(O_j)},\ \ \ R\to\infty.\tag\ehg
\endalign
$$
Note that the canonical rearrangement $\widehat{O^{(R)}}$ is the union of the individual canonical rearrangements $\widehat{O}_i$, each of which consists of the same number of monomers as separations (i.e., has charge 0). This implies
$$
E\left(\widehat{O^{(R)}}\right)\sim\prod_{i=1}^m E(\widehat{O}_i),\ \ \ R\to\infty.\tag\ehh
$$
On the other hand, the clusters in $O^{(R)}$ have in general non-zero charges, and the definition of $E$ and the hypothesis about the sequences $x_i^{(R)}$ imply
$$
E\left({O^{(R)}}\right)\sim\prod_{i=1}^m E({O}_i)
\prod_{1\leq i<j\leq m}(Rx_j-Rx_i)^{\frac12 \q(O_i)\q(O_j)},\ \ \ R\to\infty.\tag\ehi
$$
Substituting (\ehh) and (\ehi) into the left hand side of (\ehg), one readily sees that the resulting expression agrees, at least in the factors besides the powers of 2, with the one on the right hand side of (\ehg). A simple arithmetical calculation shows that the two expressions also agree in the factors that are powers of two, which establishes (\ehg) in this case. The remaining three cases ($P\geq N$, $p_u <p_v$; $P\leq N$, $p_u \geq p_v$; and $P\leq N$, $p_u <p_v$) are proved the same way. (The only differences from the above case are slight changes in the exponent of 2; compare formulas (\efa) and (\efb).)

\medskip
{\smc Step 4.} {\it Proof of the standard case when there is at least one odd cluster.} We say that an odd cluster has {\it rank $i$} if it is the $i$th cluster from the right. If a collection of defect clusters contains at least one odd cluster, the {\it rank of the collection} of defect clusters is defined to be the maximum rank of an odd cluster. When there is at least one odd cluster, we prove the standard case by induction on the rank of the collection of clusters. 

{\it Base case.} The base case is when the rank of the collection of defect clusters is 1, i.e., for each $1\leq i\leq m-1$, $\s(O_i(x_i^{(R)}))$ consists of an even number of consecutive integers, while $\s(O_m(x_m^{(R)}))$ consists of an odd number of consecutive integers. As in the previous step, we may assume, by Lemma {\tbc} and (\ebuu), that for $i=1,\dotsc,m-1$, $O_i(x_i^{(R)})$ has nice even support, and that $O_m(x_m^{(R)})$ has nice odd support. This implies that also $O_1(x_1^{(R)})\cup\cdots\cup O_m(x_m^{(R)})$ has nice odd support. Therefore, Proposition {\tfa}(b) can be applied to express all correlations occurring in (\ecb) as explicit products, and proving (\ecb) amounts to verifying that the resulting expressions on the two sides have the same $R\to\infty$ asymptotics.

As in Step 3, since our collection of defect clusters is assumed to be standardized, it follows that there is and index $u$, $1\leq u\leq m$, so that $O_1,\dotsc,O_{u-1}$ consist entirely of monomers,  $O_{u+1},\dotsc,O_{m}$ consist entirely of separations, and $O_u$ consists of a run of say $p_u$ monomers followed by a run of say $n_u$ separations, where $p_u,n_u\geq0$. In order to write down our formulas explicitly, we need to know whether $u<m$ or $u=m$ (this is due to the fact that clusters of nice even and nice odd supports have correlations given by different formulas, as stated in Proposition {\tfa}). We detail below the case $u<m$; the case $u=m$ is analogous.

Let $p_i$ be the number of monomers in $O_i$, for $i=1,\dotsc,u-1$, let $n_j$ be the number of separations in~$O_j$, for $j=u+1,\dotsc,m-1$, and denote the number of separations in the last cluster $O_m$ by $n_m+1$, so that $n_m$ is even. 

Set $O^{(R)}:=O_1(x_1^{(R)})\cup\cdots\cup O_1(x_1^{(R)})$, $P:=\sum_{i=1}^u p_i$, and $N:=\left(\sum_{j=u}^m n_j\right)+1$. Again, in order to apply Proposition {\tfa}(b), one also needs to know whether $P\geq N$ or $P\leq N$, and whether $p_u\geq n_u$ or $p_u\leq n_u$. Suppose $P\geq N$, and $p_u\geq n_u$. 

To state explicitly what (\ecb) amounts to when applying Proposition {\tfa}(b), we need some more notation. Let 
$$
\spreadlines{4\jot}
\align
S_i&:=\{s\in\Z:2s+1\in\s(O_i(x_i^{(R)}))\},\ \ \ i=1,\dotsc,m-1\tag\ehj
\\
t&:=\max\{s\in\Z:2s+1\in\s(O_m(x_m^{(R)}))\}\tag\ehk
\\
S_m&:=\{s\in\Z:s<t,2s+1\in\s(O_m(x_m^{(R)}))\}\tag\ehl
\endalign
$$

With these notations, it follows by Proposition {\tfa}(b) that in this case (\ecb) amounts to

$$
\spreadlines{4\jot}
\align
2^{         \frac{\left(\frac{P-N}{2}\right)^2}{2}           }
\prod_{i=1}^{m}\prod_{s\in S_i}\frac{\left(\frac12\right)_{t-s}^2}{(1)_{t-s}(1)_{t-s-1}}
E\left(\widehat{O^{(R)}}\right)E\left({O^{(R)}}\right)
&\sim
\frac{2^\frac{\frac{p_u-n_u}{2}\left(\frac{p_u-n_u}{2}-1\right)}{2}}{\pi^\frac{p_u+n_u}{2}}
E(\widehat{O}_u)E(O_u)
\\
&
\times
\prod_{i=1}^{u-1}\frac{2^\frac{\frac{p_i}{2}\left(\frac{p_i}{2}-1\right)}{2}}{\pi^\frac{p_i}{2}}
E(\widehat{O}_i)E(O_i)
\\
&
\times
\prod_{i=u+1}^{m-1}\frac{2^\frac{\frac{n_i}{2}\left(\frac{n_i}{2}+1\right)}{2}}{\pi^\frac{n_i}{2}}
E(\widehat{O}_i)E(O_i)
\\
&
\times
2^\frac{\frac{n_m}{2}\left(\frac{n_m}{2}+1\right)}{2}
\prod_{s\in S_m}\frac{\left(\frac12\right)_{t-s}^2}{(1)_{t-s}(1)_{t-s-1}}E(\widehat{O}_m)E(O_m)
\\
&
\times
\prod_{1\leq i<j\leq m} \left(\sqrt{2}(Rx_j-Rx_i)\right)^{\frac12 \q(O_i)\q(O_j)},\ \ \ R\to\infty.
\\
\tag\ehm
\endalign
$$
Separate the double product on the left hand side above into
$$
\prod_{i=1}^{m}\prod_{s\in S_i}\frac{\left(\frac12\right)_{t-s}^2}{(1)_{t-s}(1)_{t-s-1}}
=
\left(\prod_{i=1}^{m-1}\prod_{s\in S_i}\frac{\left(\frac12\right)_{t-s}^2}{(1)_{t-s}(1)_{t-s-1}}\right)
\prod_{s\in S_m}\frac{\left(\frac12\right)_{t-s}^2}{(1)_{t-s}(1)_{t-s-1}}.\tag\ehn
$$
Then in the limit $R\to\infty$, all the indices $t-s$ in the Pochhammer symbols involved in the factor of the double product on the right hand side in (\ehn) also tend to infinity. On the other hand, it readily follows by Stirling's formula that
$$
\lim_{r\to\infty}\frac{\left(\frac12\right)_{r}^2}{(1)_{r}(1)_{r-1}}=\frac{1}{\pi}.\tag\eho
$$
Therefore, using also the fact that $\sum_{i=1}^m |S_i|=\frac{P+N}{2}$, we obtain that
$$
\lim_{R\to\infty}\prod_{i=1}^{m}\prod_{s\in S_i}\frac{\left(\frac12\right)_{t-s}^2}{(1)_{t-s}(1)_{t-s-1}}
=\frac{1}{\pi^{\frac{P+N}{2}}}
\prod_{s\in S_m}\frac{\left(\frac12\right)_{t-s}^2}{(1)_{t-s}(1)_{t-s-1}}.\tag\ehp
$$
The asymptotics of $\bar\omega({O^{(R)}})$ follows the same way as it did in Step 3: It is
$$
E\left({O^{(R)}}\right)\sim\prod_{i=1}^m E({O}_i)
\prod_{1\leq i<j\leq m}(Rx_j-Rx_i)^{\frac12 \q(O_i)\q(O_j)},\ \ \ R\to\infty.\tag\ehq
$$
Concerning the asymptotics of $\bar\omega(\widehat{O^{(R)}})$, note the following difference from the situation in Step 3: $\widehat{O^{(R)}}$ is still the union of the individual canonical rearrangements $\widehat{O}_i$'s, but now the $\widehat{O}_i$'s have charge 0 only for $i=1,\dotsc,m-1$, while the charge of $\widehat{O}_m$ is 1. However, by the definition (\ebh) of $E$, this difference does not change the formula that gives the asymptotics of $\bar\omega(\widehat{O^{(R)}})$, since we only have one cluster of non-zero charge. Thus we have
$$
E\left(\widehat{O^{(R)}}\right)\sim\prod_{i=1}^m E(\widehat{O}_i),\ \ \ R\to\infty.\tag\ehr
$$
Substituting (\ehp)-(\ehr) into the left hand side of (\ehm), one readily sees that, except possibly for the powers of two, all factors on the two sides agree. A simple arithmetical calculation verifies that the exponents of 2 match as well, which proves (\ehm) in the case $P\geq N$, $p_u \geq p_v$. The remaining cases ($P\geq N$, $p_u <p_v$; $P\leq N$, $p_u \geq p_v$; and $P\leq N$, $p_u <p_v$) are proved the same way (again, the only differences from the above case are slight changes in the exponent of 2; compare formulas (\efc) and (\efd)).

{\it Induction step.} Suppose the collection of defect clusters $O_1,\dotsc,O_m$ is standard, and that it has rank $i>1$. Let $O_i$ be its leftmost odd cluster; thus $i<m$. Then $\s(O_1),\dotsc,\s(O_{i-1})$ consist of even-length runs of consecutive integers, and $\s(O_i)$ consists of an odd-length run of consecutive integers; $\s(O_{i+1}),\dotsc,\s(O_m)$ consist of runs of integers of even or odd length. By Lemma {\tbc} and (\ebuu), we may assume that all runs of sites separating consecutive defect clusters have even length.

The main idea of the induction step is the following. Gradually migrate the rightmost defect in $O_i$ to the right, two units at a time, until it occupies the site just to the left of the leftmost defect in $O_{i+1}$. Denote by $\tilde{O}_{i+1}$ the augmented cluster obtained this way from $O_{i+1}$, and by $\tilde{O}_i$ the cluster obtained from $O_i$ by removing its rightmost defect. Then the ratio
$$
\frac
{\bar\omega\left(O_1(x_1^{(R)}),\dotsc,O_m(x_m^{(R)})\right)}
{\bar\omega\left(O_1(x_1^{(R)}),\dotsc,O_{i-1}(x_{i-1}^{(R)}),\tilde{O}_i(x_i^{(R)}),\tilde{O}_{i+1}(x_{i+1}^{(R)}),O_{i+2}(x_{i+2}^{(R)}),\dotsc,O_m(x_m^{(R)})\right)}\tag\ehs
$$
can be expressed as an explicit product using Corollary {\tbd}. However, note that the collection of defect clusters at the denominator above is also standard, and has rank strictly less than $i$. Thus, by the induction hypothesis, the correlation at the denominator above satisfies (\ecb). Expressing the numerator in (\ehs) as the denominator times the explicit product mentioned above, checking that it satisfies (\ecb) will boil down to verifying an identity between two explicit products of ratios of linear factors.

In order to write down explicitly the resulting expressions, we need to know how the index $i$ of the leftmost odd cluster compares to the index $u$ for which clusters $O_1,\dotsc,O_{u-1}$ consist entirely of monomers, clusters $O_{u+1},\dotsc,O_m$ consist entirely of separations, and $O_u$ consists of say $p_u$ monomers and $n_u$ separations, with $p_u,n_u\geq0$ (the existence of such an index $u$ is guaranteed by the fact that we are in the standard case). 
There are a total of four cases: $i\leq u-2$ (when both $O_i$ and $O_{i+1}$ consist entirely of monomers), $i=u-1$, $i=u$, and $i>u$ (when both $O_i$ and $O_{i+1}$ consist entirely of separations). We treat in detail below the case $i=u-1$; the other cases lead only to slight changes in the resulting expressions, and are verified the same way.

\topinsert
\centerline{\mypic{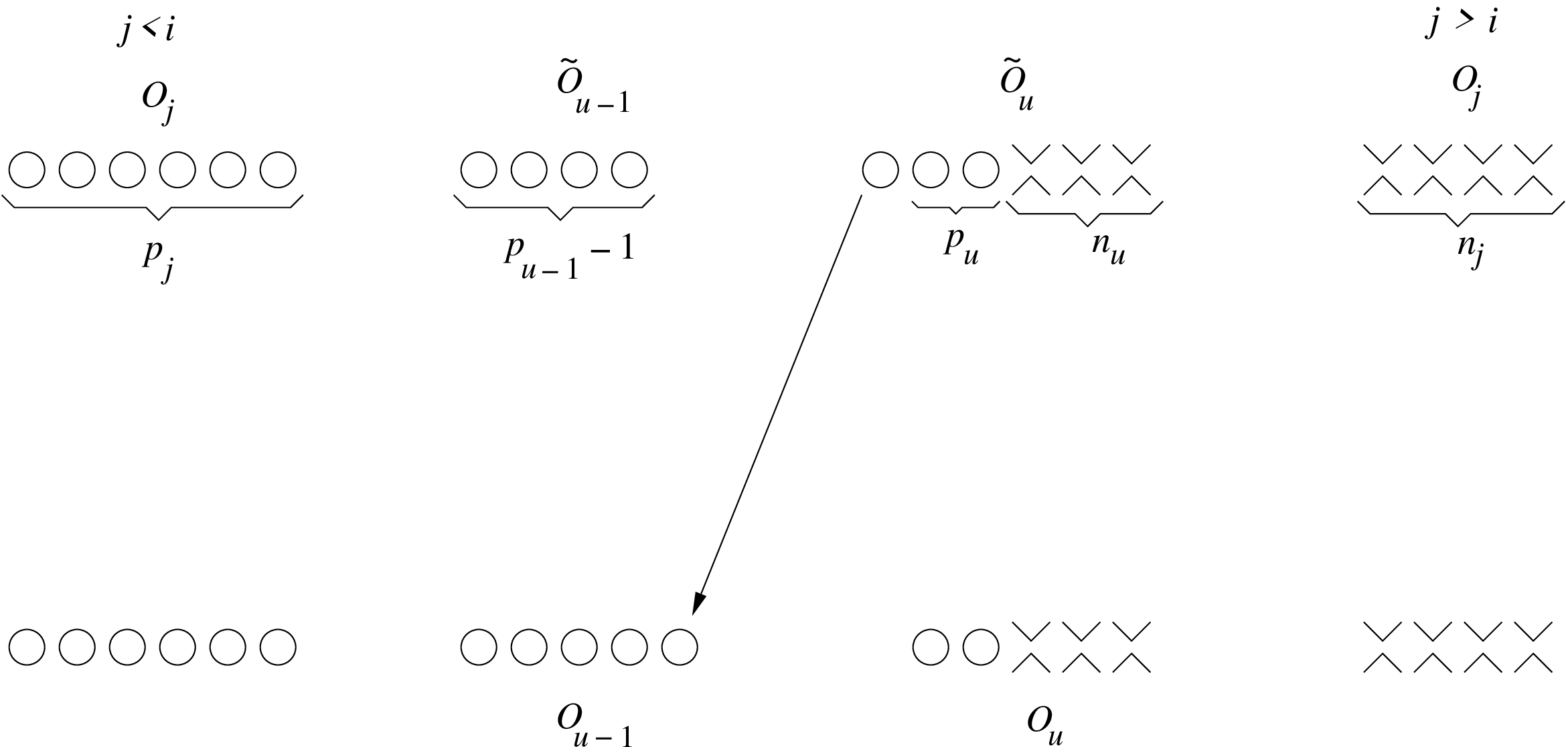}}
\medskip
\centerline{{\smc Figure~{\fha}.} Illustration of the induction step when $i=u-1$.}
\endinsert

Suppose therefore that $i=u-1$. Let $O_j$ consist of $p_j$ consecutive monomers, for $j=1,\dotsc,u-1$; then by assumption $p_1,\dotsc,p_{u-2}$ are even, and $p_{u-1}$ is odd. Let $O_j$ consist of $n_j$ consecutive separations, for $j=u+1,\dotsc,m$ (as we have mentioned above, $O_u$ consists of $p_u$ consecutive monomers followed by $n_u$ consecutive separations, with $p_u,n_u\geq0$). 

By repeated application of Corollary {\tbd}, we claim that we obtain
$$
\spreadlines{4\jot}
\align
&
\frac
{\bar\omega\left(O_1(x_1^{(R)}),\dotsc,O_{u-2}(x_{u-2}^{(R)}),\tilde{O}_{u-1}(x_{u-1}^{(R)}),\tilde{O}_{u}(x_{u}^{(R)}),\dotsc,O_m(x_m^{(R)})\right)}
{\bar\omega\left(O_1(x_1^{(R)}),\dotsc,O_m(x_m^{(R)})\right)}
\\
&
\sim
\frac
{\prod_{k=3}^{p_{u-1}+1}[L(k)]_{\lfloor\frac{Rx_u-Rx_{u-1}}{2}\rfloor}
\prod_{j=1}^{u-2}\prod_{k=3}^{p_j+2}
\dfrac{[L(k)]_{\lfloor\frac{Rx_u-Rx_{j}}{2}\rfloor}}{[L(k)]_{\lfloor\frac{Rx_{u-1}-Rx_{j}}{2}\rfloor}}}
{
\prod_{k=3}^{p_{u}+2}[L(k)]_{\lfloor\frac{Rx_u-Rx_{u-1}}{2}\rfloor}
\prod_{k=p_u+3}^{p_{u}+n_u+2}[U(k)]_{\lfloor\frac{Rx_u-Rx_{u-1}}{2}\rfloor}
\prod_{j=u+1}^{m}\prod_{k=3}^{n_j+2}
\dfrac{[U(k)]_{\lfloor\frac{Rx_j-Rx_{u-1}}{2}\rfloor}}{[U(k)]_{\lfloor\frac{Rx_{j}-Rx_{u}}{2}\rfloor}}
},\ \ \ R\to\infty\tag\eht
\endalign
$$
(recall that $[f(a)]_k$ is defined by (\efo)).

To see this, in view of the fact that Corollary {\tbd} gives the change when a defect moves two units to the left, it will be helpful to think of the collection of defect clusters at the denominator on the left hand side of (\eht) as being obtained from the collection at the numerator by migrating the leftmost monomer in $\tilde{O}_u$ gradually to the left, until it occupies the site just to the right of $\tilde{O}_{u-1}$; this causes 
$\tilde{O}_{u-1}$ to turn into $O_{u-1}$, and $\tilde{O}_{u}$ to turn into $O_{u}$.

Then the first product at the numerator on the right hand side of (\eht) arises as the contributions to the right hand side of (\ebw) coming from the monomers making up $\tilde{O}_{u-1}$ ($p_{u-1}-1$ in number), and the inner product in the double product at the numerator comes from the contributions of the monomers that constitute $O_j$, for $j=1,\dotsc,u-2$. Similarly, the first two products at the denominator on the right hand side of (\eht) come from the monomers of $\tilde{O}_u$ besides the one being moved, and from the separations of $\tilde{O}_u$, respectively, while the inner product in the double product at the denominator comes from the separations that make up $O_j$, $j=u+1,\dotsc,m$. 

The actual indices of the square brackets that follow from Corollary {\tbd} are not precisely the ones displayed in (\eht) (otherwise (\eht) would be an exact, not only an asymptotic equality). For instance, the exact value of the first index resulting by repeated application of (\ebw) is equal to half the number of sites separating $O_{u-1}$ and $O_u$. However, since $\lim_{R\to\infty}x_j^{(R)}/R=x_j$, $j=1,\dotsc,m$, one sees, using also (\ebuu), that the expression with the actual indices of the square brackets (which is {\it equal} to the left hand side of (\eht)) has the same $R\to\infty$ asymptotics as the one shown on the right in (\eht).

Applying Proposition {\tga} to the products of $L$'s and $U$'s, we obtain that the right hand side of (\eht) has asymptotics, for $R\to\infty$, given by
$$
\spreadlines{6\jot}
\align
&
\frac
{
{\!\!\!\!\!\!\!\!\!\!\!\!\!\!\!\!\!\!\!\!\!\!\!\!\!\!\!\!\!\!\!\!
p_{u-1}-1 \text{\rm \ factors} \atop \overbrace
{\displaystyle
\frac{C}{[L(3)]_0}\frac{C_e}{[L(4)]_0}\frac{C}{[L(3)]_1}\frac{C_e}{[L(4)]_1}\cdots
}
{\displaystyle
\ \left(\sqrt{\frac{Rx_u-Rx_{u-1}}{2}}\right)^{p_{u-1}-1}
}
}
}
{
{p_{u} \text{\rm \ factors} \atop \overbrace
{\displaystyle
\frac{C}{[L(3)]_0}\frac{C_e}{[L(4)]_0}\frac{C}{[L(3)]_1}\frac{C_e}{[L(4)]_1}\cdots
}
{\displaystyle
}
}
\
{\!\!\!\!\!\!\!\!\!\!\!\!\!\!\!\!\!\!\!\!\!\!\!\!\!\!\!\!\!\!\!\!n_{u} \text{\rm \ factors} \atop \overbrace
{\displaystyle
\frac{C'}{[U(3)]_{p_u}}\frac{C'_e}{[U(4)]_{p_u}}\frac{C'}{[U(3)]_{p_u+1}}\frac{C'_e}{[U(4)]_{p_u+1}}\cdots
}
{\displaystyle
\ \left(\sqrt{\frac{Rx_u-Rx_{u-1}}{2}}\right)^{p_{u}-n_u}
}
}
}
\\
&\ \ \ \ \ \ \ \ \ \ \ \ \ \ \ \ \ \ \ \ \ \ \ \ \ \ \ 
\times
\frac
{\displaystyle
\prod_{j=1}^{u-2}\left(
\frac{\sqrt{\frac{Rx_u-Rx_j}{2}}}{\sqrt{\frac{Rx_{u-1}-Rx_j}{2}}}
\right)^{p_j}
}
{\displaystyle
\prod_{j=u+1}^{m}\left(
\frac{\sqrt{\frac{Rx_j-Rx_{u-1}}{2}}}{\sqrt{\frac{Rx_{j}-Rx_u}{2}}}
\right)^{-n_j}
}.\tag\ehu
\endalign
$$
In order to write down explicitly the expression above, we need to know the parities of $p_u$ and $n_u$ (we know that $p_{u-1}-1$ is even, because $i=u-1$, and thus $p_{u-1}=p_i$ is odd). Let us assume that $p_u$ is even and $n_u$ is odd; the other parities lead only to slight changes in the resulting expressions, and are handled similarly.

By (\eht) and (\ehu) it follows then, using (\egf)-(\egh), that 
$$
\spreadlines{4\jot}
\align
&
\frac
{\bar\omega\left(O_1(x_1^{(R)}),\dotsc,O_{u-2}(x_{u-2}^{(R)}),\tilde{O}_{u-1}(x_{u-1}^{(R)}),\tilde{O}_{u}(x_{u}^{(R)}),\dotsc,O_m(x_m^{(R)})\right)}
{\bar\omega\left(O_1(x_1^{(R)}),\dotsc,O_m(x_m^{(R)})\right)}
\\
&\ \ \ \ \ \ \ \ \ 
\sim
\frac{1}{C}\,
\frac
{\displaystyle \left(\frac{4}{\pi}\right)^\frac{p_{u-1}-1}{2}}
{\displaystyle \left(\frac12\right)^\frac{n_u+1}{2}\left(\frac{4}{\pi}\right)^\frac{p_u+n_{u}-1}{2}}
\left(\sqrt{\frac{Rx_u-Rx_{u-1}}{2}}\right)^{p_{u-1}+n_u-p_u-1}
\\
&\ \ \ \ \ \ \ \ \ 
\times
\prod_{j=1}^{u-2}
\left(\sqrt{\frac{Rx_u-Rx_{j}}{2}}\right)^{p_j}
\left(\sqrt{\frac{Rx_{u-1}-Rx_{j}}{2}}\right)^{-p_j}
\\
&\ \ \ \ \ \ \ \ \ 
\times
\prod_{j=u+1}^{m}
\left(\sqrt{\frac{Rx_j-Rx_{u}}{2}}\right)^{-n_j}
\left(\sqrt{\frac{Rx_{j}-Rx_{u-1}}{2}}\right)^{n_j}
\\
&\ \ \ \ \ \ \ \ \ 
\times
\frac
{\displaystyle
\left(\prod_{j=1}^{\frac{p_u}{2}-1}[L(3)]_j[L(4)]_j\right)
\left(\prod_{j=\frac{p_u}{2}}^{\frac{p_u+n_u-1}{2}-1}[U(3)]_j[U(4)]_j\right)
[U(3)]_{\frac{p_u+n_u-1}{2}}
}
{\displaystyle
\prod_{j=1}^{\frac{p_{u-1}-1}{2}-1}[L(3)]_j[L(4)]_j
},\ \ \ R\to\infty.
\tag\ehv
\endalign
$$
On the other hand, by the induction hypothesis we know that the $R\to\infty$ asymptotics of the numerator on the left hand side of (\eht) is given by (\ecb). Besides the latter asymptotic equality, consider also the asymptotic equality we need to prove for the induction step, namely that (\ecb) holds for the {\it denominator} on the left hand side of (\eht). Taking the ratios of this two equalities, one sees, after doing the arithmetic, that (\ecb) holds for the denominator on the left hand side of (\eht) if and only if
$$
\align
\frac{\bar\omega(\tilde{O}_{u-1})\bar\omega(\tilde{O}_{u})}{\bar\omega(O_{u-1})\bar\omega(O_{u})}
\left(\sqrt{2}(Rx_u-Rx_{u-1})\right)^{\frac12(p_u-1-p_u+n_u-1)}
&\prod_{j=1}^{u-2}\left(\sqrt{2}(Rx_{u-1}-Rx_{j})\right)^{-\frac12 p_j}
                  \left(\sqrt{2}(Rx_{u}-Rx_{j})\right)^{\frac12 p_j}
\\
&\!\!\!\!\!\!\!\!\!\!\!\!\!\!
\times
\prod_{j=u+1}^{m}\left(\sqrt{2}(Rx_{j}-Rx_{u-1})\right)^{\frac12 n_j}
                  \left(\sqrt{2}(Rx_{j}-Rx_{u})\right)^{-\frac12 n_j}
\tag\ehw
\endalign
$$
has the same $R\to\infty$ asymptotics as the right hand side of (\ehv). The factors involving $Rx_l-Rx_k$ visibly agree in (\ehw) and the right hand side of (\ehv), for all $1\leq k<l\leq m$. To verify the agreement of the remaining parts of the two expressions, note that all four correlations in the fraction in (\ehw) can be expressed by the explicit product formulas provided by Proposition {\tfa}. Indeed, by our assumptions, $O_{u-1}$ and $O_u$ have nice odd support, while $\tilde{O}_{u-1}$ and $\tilde{O}_u$ have nice even support (note that repeated application of the elementary move lemma implies that translating a defect cluster one unit to the left keeps its correlation unchanged). 

Therefore, using formulas (\efa)-(\efd) to express the correlations in (\ehw) (note that, depending on whether $p_u\geq n_u$ or $p_u<n_u$, we need to use (\efc) or (\efd) for $\bar\omega(O_u)$, and (\efa) or (\efb) for $\bar\omega(\tilde{O}_u)$, respectively), checking the agreement of the remaining parts of (\ehw) and the right hand side of (\ehv) amounts to verifying the equality of two explicit products of ratios of Gamma functions. Using also (\egk) and (\egl) to express the products of $U$'s in terms of products of $L$'s, it is routine to check that the remaining parts of (\ehw) and the right hand side of (\ehv) agree. This completes the proof in the case when the leftmost odd cluster of the collection $O_1,\dotsc,O_m$ is $O_{u-1}$.

The remaining three cases (leftmost odd cluster is $O_i$ with $i<u-1$, $i=u$, resp. $i>u$) entail only slight changes to the calculations above, and are verified in the same way. \epf

\medskip
\flushpar
{\smc Remark 9.} It is interesting to note that a weaker version of Theorem {\tca}, in which the multiplicative constant on the right hand side of (\eca) is not specified, can be obtained more directly using only exactness, the elementary move lemma, and the asymptotics of the products of $L$'s and $U$'s. Indeed, the standard case can be then proved by gradually migrating all defect clusters until they coalesce, keeping track of the change in correlation by using the elementary move lemma. The exact formulas for the correlation given in Section 6 and the inductive argument above are needed however in order to identify the multiplicative constant as the product of the correlations of the individual clusters.

\mysec{9. Concluding remarks}

The main result of this paper, Theorem {\tca}, proves a version of the conjectured strong superposition principle~(\ebe) in the case when the defect clusters are collinear. This is the first proved evidence that involves defects of arbitrary charge (even or odd) for (\ebe), whose analog on the hexagonal lattice was conjectured in \cite{\ov,Conjecture\ 1}.

Our proof consists of the asymptotic analysis of an exact product formula (see Theorem {\taa}) that gives the number of perfect matchings of Aztec rectangles with defects of two kinds~---~monomers and separations~---~along a symmetry axis. 

A seemingly unrelated but well-known problem is the question of determining the interaction of two bubbles in a universe filled with a homogeneous fluid, any two points of which attract according to Newtonian gravity. The answer is that the bubbles also attract. What we have found in the current paper is an analog of this phenomenon for electrostatics.

Indeed, briefly speaking, we have found that if we perturb the regular sequence of alternating $O$'s and $E$'s on $\Z$ by introducing unit gaps corresponding to monomers, and unit overlaps corresponding to separations, and we impose two independent, repulsive Coulomb interactions, one on the sites occupied by the $O$'s, the other on the sites occupied by the $E$'s, then the resulting effect in the scaling limit is a Coulomb interaction between the {\it defects}, in which monomers act as unit charges of one sign, and separations as unit charges of the opposite sign. Moreover, the arising Coulomb interaction has exponent $\frac12$, in the sense that for instance the correlation of a monomer and a separation distance $d$ apart goes to zero as $d^{-\frac12}$, when $d\to\infty$.

The simpler (but not subsumed by the above) case of having just $O$'s at all locations in $\Z$, perturbing this sequence by introducing unit gaps, and determining the effect on the gaps of having a repulsive Coulomb interaction on the $O$'s, was worked out in equivalent language in \cite{\sc}. The result we found there was that the gaps were in their turn governed by a repulsive Coulomb interaction, of the same strength as the original one.

It is natural to consider the same problem also for other patterns, say perturbations of $\Z^2$ or $\Z^3$ filled with $O$'s caused by introducing gaps, or perturbations of more general patterns containing more than one symbol. Even though currently there seem to be no combinatorial set-ups leading to such patters, it would still be interesting to analyze them and see what kind of interactions on the defects they lead to.

%
%
%
%
%
%
%
%
%
%
%
%
%
%
%
%
%
%
%
%
%
%
%
%
%
%

\bigskip

\mysec{References}
{\openup 1\jot \frenchspacing\raggedbottom
\roster

\myref{\Baikone}
  J. Baik, T. Kriecherbauer, K. McLaughlin and P. Miller, Uniform asymptotics for polynomials orthogonal with respect to a general class of discrete weights and universality results for associated ensembles: Announcement of results, {\it Int. Math. Res. Not.} {\bf 2003} (2003), 821--858.

\myref{\Baiktwo}
  J. Baik, T. Kriecherbauer, K. McLaughlin and P. Miller, ``Discrete orthogonal polynomials. Asymptotics and applications,'' {\it Annals of Math Studies,} Princeton University Press, Princeton, NJ, 2007. 

\myref{\FT}
  M. Ciucu, Enumeration of perfect matchings in graphs with reflective
symmetry, {\it J. Comb. Theory Ser. A} {\bf 77} (1997), 67--97.

\myref{\ri}
  M. Ciucu, Rotational invariance of quadromer correlations on the hexagonal lattice, {\it Adv. in Math.} {\bf 191} 
(2005), 46-77.

\myref{\sc}
  M. Ciucu, A random tiling model for two dimensional electrostatics, {\it Mem. Amer. Math. Soc.} {\bf 178} (2005),
no. 839, 1--106.

\myref{\ppone}
  M. Ciucu, A random tiling model for two dimensional electrostatics, {\it Mem. Amer. Math. Soc.} {\bf 178} (2005),
no. 839, 107--144.

\myref{\ec}
  M. Ciucu, The scaling limit of the correlation of holes on the triangular lattice with periodic boundary 
conditions, {\it Mem. Amer. Math. Soc.} {\bf 199} (2009), 1--100.

\myref{\ef}
  M. Ciucu, The emergence of the electrostatic field as a Feynman sum in random tilings with holes, {\it  Trans. Amer. Math. Soc.} {\bf 362}  (2010), 4921--4954.

\myref{\ov}
  M. Ciucu, Dimer packings with gaps and electrostatics, {\it Proc. Natl. Acad. Sci. USA}  {\bf 105}  (2008),  no. 8, 2766--2772.

\myref{\CEP}
  H. Cohn, N. Elkies, and J. Propp, Local statistics for random domino tilings of the 
Aztec diamond, {\it Duke Math. J.} {\bf 85} (1996), 117-166.

%
%

\myref{\FS} 
  M. E. Fisher and J. Stephenson, Statistical mechanics of dimers on a plane 
lattice. II. Dimer correlations and monomers, {\it Phys. Rev. (2)} {\bf 132} (1963),
1411--1431.

\myref{\FisherIsing}
  M. E. Fisher, On the dimer solution of planar Ising models, {\it J. Math. Phys.} {\bf 7} (1966), 1776--1781.

\myref{\Glaish}
  J. W. L. Glaisher, On Certain Numerical Products in which the Exponents Depend Upon the Numbers, {\it  Messenger Math.} {\bf 23} (1893), 145--175.


\myref{\Hart} 
  R. E. Hartwig, Monomer pair correlations, {\it J. Mathematical Phys.} {\bf 7}
(1966), 286--299.

\myref{\SM}
  W. Krauth and  R. Moessner, Pocket Monte Carlo algorithm for classical doped dimer models,
{\it Physical Review B} {\bf 67} (2003), 064503.

\myref{\MRR}
  W. H. Mills, D. P. Robbins, and H. Rumsey, Alternating sign matrices and
descending plane partitions, {\it J. Comb. Theory Ser. A} {\bf 34} (1983), 
340--359.

\myref{\ZI}
  J. B. Zuber and C. Itzykson, Quantum field theory and the two-dimensional Ising model, {\it Phys. Rev. D} {\bf 15} (1977), 2875--2884.

\endroster\par}

\enddocument